# COVARIATE ADJUSTED FUNCTIONAL PRINCIPAL COMPONENTS ANALYSIS FOR LONGITUDINAL DATA

By Ci-Ren Jiang[1] and Jane-Ling Wang[2]

*University of California, Davis*

Classical multivariate principal component analysis has been extended to functional data and termed functional principal component analysis (FPCA). Most existing FPCA approaches do not accommodate covariate information, and it is the goal of this paper to develop two methods that do. In the first approach, both the mean and covariance functions depend on the covariate $Z$ and time scale $t$ while in the second approach only the mean function depends on the covariate $Z$. Both new approaches accommodate additional measurement errors and functional data sampled at regular time grids as well as sparse longitudinal data sampled at irregular time grids. The first approach to fully adjust both the mean and covariance functions adapts more to the data but is computationally more intensive than the approach to adjust the covariate effects on the mean function only. We develop general asymptotic theory for both approaches and compare their performance numerically through simulation studies and a data set.

**1. Introduction.** Principal component analysis is a standard dimension reduction tool for multivariate data and has been extended to functional data that are in the form of random curves. Because functional data are intrinsically infinite dimensional, dimension reduction is essential to analyze such data. In addition to Ferraty and Vieu (2006) and Wu and Zhang (2006), the sequence of monographs by Ramsay and Silverman (2002, 2005) provide a tutorial for the methodology and applications of "Functional Data Analysis" (FDA). A sizable literature exists for FPCA, functional approaches to conduct principal component analysis, when entire curves are observed for each subject or in practical terms when subjects are measured at

Received January 2009; revised July 2009.

[1]Based in part on Ph.D. dissertation supported by NSF Grant DMS-04-06430.

[2]Supported in part by NSF Grant DMS-04-06430 and NIH Grants P01-DK-45939 and P01-AG022500-01.

a dense grid of time points [see, e.g., Rao (1958), Dauxois, Pousse and Romain (1982), Besse and Ramsay (1986), Castro, Lawton and Sylvestre (1986), Rice and Silverman (1991), Boente and Fraiman (2000), Bosq (2000), Cardot (2000, 2006), Mas and Menneteau (2003) and Hall and Hosseini-Nasab (2006)]. Kneip and Utikal (2001) used methods of FDA to assess the variability of densities for data sets from different populations. When functional data are observed at irregular time points, perhaps just a few time points per subject, they are usually referred as longitudinal data since they often arise from longitudinal studies. Rice (2004) and Hall, Müller and Wang (2006) described the intrinsic similarities and differences between FDA and longitudinal data analysis. Longitudinal data are often sparse with few measurements per subject and noisy with measurement errors (or random fluctuations). However, these difficulties can be overcome in most situations, so it is still possible to conduct FPCA [Shi, Weiss and Taylor (1996), James, Hastie and Suger (2000), Rice and Wu (2001), Yao, Müller and Wang (2005), Paul and Peng (2009) and Peng and Paul (2009)].

The aforementioned FPCA approaches treat all functional data as if they come from the same population. When the covariate information is available, some non-FPCA approaches such as functional mixed effects models [Wang (1998) and Guo (2002)] and semiparametric mixed effects models Zhang et al. (1998) are proposed. There has been little work involving covariate information in the framework of FPCA although it might be of particular interest in many situations, for example, to study the modes of variation of the data. Furthermore, FPCA is an effective dimension reduction method. Chiou, Müller and Wang (2003) considered a general approach incorporating a vector covariate effect through a semiparametric model. Their approach consists of two steps. In the first step, traditional FPCA was performed on all subjects ignoring the covariate information. This resulted in a Karhunen–Loève expansion [see (2.2)] for each subject $X(t)$ for which the conditional expectation of $X(t)$ given the covariate $Z$ was obtained and subsequently estimated through a semiparametric approach. A different approach was proposed in Cardot (2006), who considered conditional FPCA through nonparametric kernel estimators of the conditional mean functions and conditional variance functions. A key assumption for both approaches is that the trajectories of the functional data are either completely observed or densely recorded over time. Both assumptions are rarely satisfied in longitudinal medical or social studies. Specifically, the approach in Chiou, Müller and Wang (2003) is not suitable for extension to sparse longitudinal data as the conditional principal components cannot be estimated or approximated consistently for sparse longitudinal data. We propose a unified approach in Section 2 to model the mean function and two different approaches to model the covariance function.

Little is known on how to incorporate covariate information in FPCA for sparse longitudinal data, so our goal in this paper is to provide a unified platform to incorporate the covariate information that is applicable to both functional and longitudinal data. Two different approaches are proposed; one is based on conditional FPCA and the other adjusts the covariate effect on the mean function only. We derive uniform consistency and asymptotical normality for the mean and covariance functions for kernel and local polynomial smoothers. The two approaches are compared numerically through a simulation study and illustrated with a data example.

The rest of this paper is organized as follows: Section 2 introduces the two new approaches and their estimation procedures. Asymptotic results and the theoretical properties of the proposed estimators are described in Section 3 with proofs in the Appendix. Practical implementations of the new approaches and simulation studies are discussed in Section 4. In Section 5, we employ both approaches to the Mexican Flies data in Carey et al. (2005) and compare them to three FPCA approaches James, Hastie and Suger (2000), Yao, Müller and Wang (2005), Peng and Paul (2009) that do not incorporate covariate information. Conclusions are in Section 6.

**2. Methodology.** Ignoring the covariate information for the moment and consider the random functions $X(t)$ with mean $\mu(t)$ and covariance $\Gamma(t,s)$. FPCA in this simple setting corresponds to a spectral decomposition of the covariance $\Gamma$ and leads to the Karhunen–Loève decomposition of the random function

$$(2.1) \qquad X(t) = \mu(t) + \sum_{k=1} A_k \phi_k(t),$$

where $\phi_k(t)$ is the eigenfunction of the covariance function $\Gamma(s,t)$ corresponding to the $k$th largest eigenvalues, and $A_k = \int_{\mathcal{T}} \{X(t) - \mu(t)\} \phi_k(t)\, dt$ is the $k$th functional principal component score. In the presence of a covariate $Z = z$ we view $X(t,z)$ as a random function with mean function $\mu(t,z)$ and covariance function $\Gamma(t,s,z)$ where $s$ and $t$ are in a compact time interval $\mathcal{T}$. In this paper, the random function $X(t,z)$ are not observable because measurements are taken on discrete time points and there may be measurement errors. This is different from the situation considered in Cardot (2006) where a covariate adjusted FPCA was proposed under the assumption that the entire function $X(t,z)$ can be observed without errors.

2.1. *Model.* We consider two ways to extend the FPCA approach to accommodate covariate information. Both approaches consist of two parts: a systematic part corresponding to the mean function and a stochastic part comprising the random components that reflect the covariance structure of the longitudinal data. In both approaches we do not assume that we know the structure of $\mu(t, z)$ other than that it is a smooth function, so we will need to estimate it nonparametrically. The difference between the two approaches is in the handling of the covariance structure. Conceptually, the covariate $Z$ can be any vector that has continuous distribution, but due to the curse of dimensionality only low-dimensional $Z$ can be used. Some dimension reduction approaches will be necessary for high-dimensional $Z$ and are beyond the scope of this paper.

In the first approach, it is assumed that the eigenfunctions of $\Gamma(t, s, z)$ vary with $z$ so that there exists an orthogonal expansion of $\Gamma$ (in the $L^2$ sense) in terms of eigenfunctions $\phi_k(t, z)$ and nonincreasing eigenvalues $\lambda_k(z): \Gamma(t, s, z) = \sum_k \lambda_k(z)\phi_k(t, z)\phi_k(s, z)$. Thus the random trajectory $X(t, z)$ can be represented as

$$(2.2) \qquad X(t, z) = \mu(t, z) + \sum_{k=1} A_k(z)\phi_k(t, z),$$

where $A_k(z)$ are uncorrelated random variables with mean 0 and variance $\lambda_k(z)$. Again, we will model the covariance surface nonparametrically, assuming that it is a smooth function of $t, s$ and $z$. Since both the mean and covariance functions have been adjusted by the covariate $Z$, we call this fully adjusted functional principal component analysis and abbreviate it as fF-PCA. This approach to adjust the covariate effects is conceptually equivalent to the conditional FPCA approach in Cardot (2006) but differs crucially in the way of estimation due to differences in the data design we consider here. This crucial difference in the data design also triggers a much different theoretical framework. For one-dimensional $Z$, only one-dimensional smoothing is needed in Cardot (2006) to estimate both the mean and covariance function along the $Z$-direction at each time location since the entire function $X(t, z)$ is observed.

When $\mu(t, z) = \beta(t)z$ and the stochastic components $\sum_{k=1} A_k(z)\phi_k(t, z)$ in model (2.2) adopts a time-varying linear structure $b(t)z$ for some unknown function $\beta$ and random function $b$, model (2.2) yields the varying coefficient random effects model in Guo (2002). When $\mu(t, z)$ takes a partial linear form $f(t) + \beta z$ and the stochastic component also takes a partial linear form $u(t) + bZ$, for some unknown functions $f$ and $u$, parameter $\beta$ and random variable $b$, model (2.2) reduces to the partial linear mixed model in Zhang et al. (1998).

In the second approach, one can take advantage of the fact that the covariate $Z$ is a random variable and pool all subjects together after centering each individual curve at zero. This leads to a pooled covariance function $\Gamma^*(t, s) = \int_{\mathcal{Z}} \mathrm{E}\{(X(t, z) - \mu(t, z))(X(s, z) - \mu(s, z))\}g(z)\,dz$ where $g$ is the p.d.f. of $Z$ on $\mathcal{Z}$, and $\Gamma^*(t, s)$ is assumed to be a smooth function of $t$ and $s$. Consequently, there exists an orthogonal expansion (in the $L^2$ sense) in terms of eigenfunctions $\phi_k^*$ and nonincreasing eigenvalues $\lambda_k^*$ such that $\Gamma^*(t,s) = \sum_k \lambda_k^* \phi_k^*(t)\phi_k^*(s)$, and

$$(2.3) \qquad X(t,z) = \mu(t,z) + \sum_{k=1} A_k^* \phi_k^*(t),$$

where $A_k^*$ are uncorrelated random variable with $E\{A_k^*\} = 0$ and $\text{var}\{A_k^*\} = \lambda_k^*$. This approach has the advantage that the covariance function can be estimated with a lower-dimensional smoother than its counterpart in fFPCA, accelerating the rate of convergence compared to fFPCA. We abbreviate this mean adjusted functional principal component analysis on $X(t,z) - \mu(t,z)$ as "mFPCA" where "m" stands for the mean adjusting operation. The estimation procedure for mFPCA is described in Section 2.2.2. Conceptually, the fFPCA approach should fit the data better as it adapts to the covariate information in covariance estimation while mFPCA does not. This benefit may be offset by inferior practical performance if the data are sparse. Our simulation results in Section 4 reflect limited benefits of fFPCA, so one may prefer the mFPCA approach in many applications or try both approaches, unless the eigenfunctions vary substantially across the covariate values.

2.2. *Estimation.* In many situations one can only observe the processes $X(t,z)$ intermittently and possibility with measurement errors. Let $Y_{ij}$ be the $j$th observation of the random function $X_i$, made at a random time $T_{ij} \in \mathcal{T}$ with a covariate $Z_i \in \mathcal{Z}$ and measurement error $\epsilon_{ij}$ where $i = 1, \ldots, n$, and $j = 1, \ldots, N_i$. Here we assume that the measurement schedule $T_{ij}$ is a random sample of size $N_i$ and $N_i$ is assumed to be i.i.d. and independent of all other random variables. We also assume that the measurement errors are i.i.d. with mean 0 and a constant variance $\sigma^2$ and are independent of the random coefficients $A_k(z)$ or $A_k^*$ under model (2.2) or (2.3), respectively. Thus the observed data are

$$(2.4) \qquad Y_{ij} = X_i(T_{ij}, Z_i) + \epsilon_{ij}.$$

The key steps in our FPCA approach are to estimate the mean and covariance function. The corresponding eigenvalues and eigenfunctions can be obtained easily through the eigen-equation after the covariance function has been estimated. The mean functions for fFPCA and mFPCA are the same and can be estimated using any two-dimensional scatter-plot smoother of $Y_{ij}$ against $(T_{ij}, Z_i)$, for $j = 1, \ldots, N_i, i = 1, \ldots, n$. We provide general asymptotic properties of any linear scatter-plot smoother of the mean function $\mu(t,z)$ and demonstrate in Section 3 these asymptotic properties on two linear smoothers, the Nadaraya–Watson estimator (3.1) and the local linear estimator (3.2).

Similarly, our covariance estimator can also be expressed as a scatter-plot smoother of the so called "raw covariances" defined below against $(T_{ij}, T_{ik})$:

$$(2.5) \qquad C_{ijk} = (Y_{ij} - \hat{\mu}(T_{ij}, Z_i))(Y_{ik} - \hat{\mu}(T_{ik}, Z_i)).$$

The covariance estimators are different for fFPCA and mFPCA. For one-dimensional $Z$, the former involves a three-dimensional smoother of $C_{ijk}$ against $(T_{ij}, T_{ik}, Z_i)$ for $j, k = 1, \ldots, N_i, i = 1, \ldots, n$ while the latter only requires a two-dimensional smoother of $C_{ijk}$ against $(T_{ij}, T_{ik})$ for $j, k = 1, \ldots, N_i, i = 1, \ldots, n$. In principle, one can employ any linear smoother. We illustrate the theorems for the Nadaraya–Watson estimators and local linear estimators [Fan and Gijbels (1996)] in Section 3.

### 2.2.1. fFPCA. Note that since

$$\operatorname{cov}(Y_{ij}, Y_{ik} | T_{ij}, T_{ik}, Z_i)$$
$$= \operatorname{cov}(X(T_{ij}, Z_i), X(T_{ik}, Z_i)) + \sigma^2 \delta_{jk},$$

where $\delta_{jk}$ is 1 if $j = k$ and 0 otherwise, the diagonal of the "raw" covariances $C_{ijk}$ in (2.5) should not be included in the covariance function smoothing step. With this in mind the local linear smoother for the covariance function $\Gamma(t, s, z)$ is

$$\hat{\Gamma}_L(t, s, z) = \hat{\beta}_0 \qquad \text{where}$$

$$
\begin{aligned}
(2.6) \qquad \hat{\boldsymbol{\beta}} = \underset{\boldsymbol{\beta}}{\arg\min} \Bigg\{ & \sum_{i=1}^{n} \sum_{1 \leq j \neq k \leq N_i} K_3 \left( \frac{t - T_{ij}}{h_{G,t}}, \frac{s - T_{ik}}{h_{G,t}}, \frac{z - Z_i}{h_{G,z}} \right) \\
& \times [C_{ijk} - (\beta_0 + \beta_1(T_{ij} - t) \\
& \qquad + \beta_2(T_{ik} - s) + \beta_3(Z_i - z))]^2 \Bigg\},
\end{aligned}
$$

and $K_3$ is a three-dimensional kernel function satisfying (A.2).

Next we aim to estimate the variance $V(t, z) = \Gamma(t, t, z) + \sigma^2$ of $Y(t)$ for a given $z$. Let $K_2$ be a two-dimensional kernel function satisfying (A.1) and $\hat{V}(t, z)$ be the local linear smoother using only the diagonal time elements; then

$$\hat{V}(t, z) = \hat{\beta}_0,$$

where

$$
\begin{aligned}
\hat{\boldsymbol{\beta}} = \underset{\boldsymbol{\beta}}{\arg\min} \sum_{i=1}^{n} \sum_{j=1}^{N_i} K_2 & \left( \frac{t - T_{ij}}{h_{V,t}}, \frac{z - Z_i}{h_{V,z}} \right) \\
& \times [C_{ijj} - \beta_0 - \beta_1(T_{ij} - t) - \beta_2(Z_i - z)]^2.
\end{aligned}
$$

The variance $\sigma^2$ of the measurement error can be estimated by averaging $(\hat{V}(t,z) - \hat{\Gamma}_L(t,t,z))$ over a range of $t$. For stability, one may prefer to use a trimmed mean restricting the averaging to take place over a central portion of the time domain. We find the recommendation in Yao, Müller and Wang (2005) to use a trimmed mean based on the central 50% of the time domain satisfactory. Specifically, this leads to

$$(2.7) \qquad \hat{\sigma}^2 = \frac{1}{|\mathcal{T}_1||\mathcal{Z}|} \int_{\mathcal{Z}} \int_{\mathcal{T}_1} \{\hat{V}(t,z) - \hat{\Gamma}_L(t,t,z)\} \, dt \, dz,$$

where $\mathcal{T}_1$ is the interval $[\inf\{t : t \in \mathcal{T}\} + |\mathcal{T}|/4, \sup\{t : t \in \mathcal{T}\} - |\mathcal{T}|/4]$ with the notation $|\mathcal{I}|$ denoting the length of a generic interval $\mathcal{I}$. If the variances of the measurement errors vary over time and $z$, the variance function $\sigma^2(t,z)$ can be estimated directly as $\hat{V}(t,z) - \hat{\Gamma}(t,t,z)$.

The solutions of the eigen-equations, $\int \hat{\Gamma}_L(t,s,z) \hat{\phi}_k(s,z) \, ds = \hat{\lambda}_k(z) \hat{\phi}_k(t,z)$, where the $\hat{\phi}_k(t,z)$ satisfies $\int \hat{\phi}_k^2(t,z) \, dt = 1$ and $\int \hat{\phi}_k(t,z) \hat{\phi}_m(t,z) \, dt = 0$ for $m < k$, are used to estimate the eigenfunctions and eigenvalues. It now remains to estimate the principal component scores

$$A_{ik}(Z_i) = \int \phi_k(t, Z_i)[X_i(t, Z_i) - \mu(t, Z_i)] \, dt$$

for the $i$th subject. Due to measurement errors and the intermittent measurement schedules, the approaches in Chiou, Müller and Wang (2003) and Cardot (2006) are not applicable to estimate these scores. Instead, the approach in Yao, Müller and Wang (2005) aimed at estimating the conditional expectation $E(A_{ik}(Z_i)|\tilde{\mathbf{Y}}_i)$ is well suited to estimate the principal component scores where $\tilde{\mathbf{Y}}_i = (Y_{i1}, \ldots, Y_{iN_i})^T$. Under the assumption that $\tilde{\mathbf{Y}}_i$ is multivariate normal, this leads to the estimate

$$\hat{A}_{ik}(Z_i) = \hat{\lambda}_k \hat{\boldsymbol{\phi}}_{ik}^T \hat{\Sigma}_{\tilde{\mathbf{Y}}_i}^{-1} (\tilde{\mathbf{Y}}_i - \hat{\boldsymbol{\mu}}_i),$$

where $\hat{\boldsymbol{\mu}}_i = (\hat{\mu}(T_{i1}, Z_i), \ldots, \hat{\mu}(T_{iN_i}, Z_i))^T$, $(\hat{\Sigma}_{\tilde{\mathbf{Y}}_i})_{j,k} = \hat{\Gamma}_L(T_{ij}, T_{ik}, Z_i) + \hat{\sigma}^2 \delta_{jk}$ and $\hat{\boldsymbol{\phi}}_{ik} = (\hat{\phi}_k(T_{i1}, Z_i), \ldots, \hat{\phi}_k(T_{iN_i}, Z_i))^T$.

2.2.2. *mFPCA*. The estimation of $\Gamma^*(s,t)$ is similar to the procedure in Yao, Müller and Wang (2005) except that we use $C_{ijk}$ as the raw covariances. Let $\hat{\Gamma}^*(t,s)$ be the covariance estimator based on a local linear smoother; then

$$\hat{\Gamma}^*(t,s) = \hat{\beta}_0 \qquad \text{where}$$

$$\hat{\boldsymbol{\beta}} = \arg\min_{\boldsymbol{\beta}} \left\{ \sum_{i=1}^{n} \sum_{1 \le j \ne k \le N_i} K_2 \left( \frac{t - T_{ij}}{h_{G^*}}, \frac{s - T_{ik}}{h_{G^*}} \right) \right.$$

(2.8)
$$\times [C_{ijk} - (\beta_0 + \beta_1(T_{ij} - t)$$
$$+ \beta_2(T_{ik} - s))]^2 \bigg\},$$

where $t, s \in \mathcal{T}$ and $K_2$ is defined in (A.1). Let $\hat{V}^*(t)$ be the local linear smoother focusing on the diagonal values $\{\Gamma^*(t,t) + \sigma^2\}$; then

$$\hat{V}^*(t) = \hat{\beta}_0,$$

where

$$\hat{\boldsymbol{\beta}} = \underset{\boldsymbol{\beta}}{\arg\min} \sum_{i=1}^{n} \sum_{j=1}^{N_i} K_1 \left( \frac{t - T_{ij}}{h_{V^*}} \right) [C_{ijj} - \beta_0 - \beta_1(T_{ij} - t)]^2,$$

where $K_1$ is a kernel function with compact support, symmetric and Hölder continuous. Again, a "trimmed" mean of $(\hat{V}^*(t) - \hat{\Gamma}^*(t,t))$ is used to estimate $\sigma^2$ similar to (2.7).

The solutions of the eigen-equations, $\int \hat{\Gamma}^*(t,s)\hat{\phi}_k^*(s)\,ds = \hat{\lambda}_k^*\hat{\phi}_k^*(t)$, where the $\hat{\phi}_k^*(t)$ satisfies $\int (\hat{\phi}_k^*(t))^2\,dt = 1$ and $\int \hat{\phi}_k^*(t)\hat{\phi}_m^*(t)\,dt = 0$ for $m < k$, are used to estimate the eigenfunctions and eigenvalues. The principal component scores $A_{ik}^*$ for subject $i$ are estimated as in Yao, Müller and Wang (2005) through

$$\hat{A}_{ik}^* = \hat{\lambda}_k^*(\hat{\boldsymbol{\phi}}_{ik}^*)^T (\hat{\Sigma}_{\tilde{\mathbf{Y}}_i}^*)^{-1}(\tilde{\mathbf{Y}}_i - \hat{\boldsymbol{\mu}}_i),$$

where $\tilde{\mathbf{Y}}_i$ and $\hat{\boldsymbol{\mu}}_i$ are defined as in Section 2.2.1, and $(\hat{\Sigma}_{\tilde{\mathbf{Y}}_i}^*)_{j,k}$ and $\hat{\boldsymbol{\phi}}_{ik}^*$ are defined as $(\hat{\Sigma}_{\tilde{\mathbf{Y}}_i}^*)_{j,k} = \hat{\Gamma}^*(T_{ij}, T_{ik}) + (\hat{\sigma}^*)^2\delta_{jk}$ and $\hat{\boldsymbol{\phi}}_{ik}^* = (\hat{\phi}_k^*(T_{i1}), \dots, \hat{\phi}_k^*(T_{iN_i}))^T$.

2.3. *Bandwidth selection and number of eigenfunctions.* The bandwidths for the estimated mean function are chosen via the leave-one-curve-out cross-validation suggested by Rice and Silverman (1991). However, the bandwidths of the covariance function estimators are chosen via a $k$-fold cross-validation procedure to save computing time. Below we define the $k$-fold cross-validation method for the bandwidths selection of $\Gamma^*(t,s)$. The formula for $\Gamma(t,s,z)$ is similar.

Suppose that the subjects are randomly assigned to $k$ sets $(S_1, S_2, \dots, S_k)$.

(2.9)    $$h = \underset{h}{\arg\min} \sum_{\ell=1}^{k} \sum_{i \in S_\ell} \sum_{1 \leq j \neq m \leq N_i} \{C_{ijm} - \hat{\Gamma}^{*(-S_\ell)}(T_{ij}, T_{im})\}^2,$$

where $\hat{\Gamma}^{*(-S_\ell)}(T_{ij}, T_{im})$ is the estimated covariance function at $(T_{ij}, T_{im})$ when the subjects in $S_\ell$ are not used to estimate $\Gamma^*(t,s)$. We found the ten-fold $(k = 10)$ method to have satisfactory performance.

Three criteria to choose the number of eigenfunctions $K$ are discussed in the simulation study section. Suppose that the first $K$ eigenfunctions are used to predict the trajectories; given $t \in \mathcal{T}$ and $z \in \mathcal{Z}$, the predicted trajectory of $X_i(t, z)$ based on the first $K$ eigenfunctions will be

$$\text{(fFPCA)} \qquad \hat{X}_i^K(t, z) = \hat{\mu}_L(t, z) + \sum_{k=1}^{K} \hat{A}_{ik}(z) \hat{\phi}_k(t, z),$$

$$\text{(mFPCA)} \qquad \tilde{X}_i^K(t, z) = \hat{\mu}_L(t, z) + \sum_{k=1}^{K} \hat{A}_{ik}^* \hat{\phi}_k^*(t).$$

## 3. Asymptotic results.
For simplicity, the covariate $Z$ in this section will be univariate, and $N_1, \ldots, N_n$ are i.i.d. copies of some random variable $N$. We first focus on the asymptotic distribution of linear smoothers of the mean function.

*Asymptotic results for mean functions.* A general theory (Lemma C.2) for two-dimensional kernel-weighted estimators is provided in the Appendix; from there the asymptotic normality of the Nadaraya–Watson kernel estimator $\hat{\mu}_{\mathrm{NW}}(t, z)$ and local linear estimator $\hat{\mu}_L(t, z)$ of $\mu(t, z)$ follows. Specifically,

$$(3.1) \quad \hat{\mu}_{\mathrm{NW}}(t, z) = \frac{\sum_{i=1}^{n} \sum_{j=1}^{N_i} K_2((t - T_{ij})/h_{\mu,t}, (z - Z_i)/h_{\mu,z}) Y_{ij}}{\sum_{i=1}^{n} \sum_{j=1}^{N_i} K_2((t - T_{ij})/h_{\mu,t}, (z - Z_i)/h_{\mu,z})},$$

$$\hat{\mu}_L(t, z) = \hat{\beta}_0 \qquad \text{where for } \boldsymbol{\beta} = (\beta_0, \beta_1, \beta_2),$$

$$(3.2) \qquad \hat{\boldsymbol{\beta}} = \operatorname*{arg\,min}_{\boldsymbol{\beta}} \sum_{i=1}^{n} \sum_{j=1}^{N_i} K_2\left(\frac{t - T_{ij}}{h_{\mu,t}}, \frac{z - Z_i}{h_{\mu,z}}\right)$$
$$\times [Y_{ij} - \beta_0 - \beta_1(T_{ij} - t) - \beta_2(Z_i - z)]^2.$$

THEOREM 3.1. *Under assumptions* A.3, A.5 *and* A.6, B.1–B.4, *and assuming* $\frac{h_{\mu,z}}{h_{\mu,t}} \to \rho_\mu$ *and* $nE(N)h_{\mu,t}^6 \to \tau_\mu^2$ *for some* $0 < \rho_\mu, \tau_\mu < \infty$,

$$\sqrt{n\bar{N}h_{\mu,t}h_{\mu,z}}[\hat{\mu}_{\mathrm{NW}}(t, z) - \mu(t, z)] \xrightarrow{\mathcal{D}} N(\beta_{\mathrm{NW}}, \Sigma_{\mathrm{NW}}),$$

*where*

$$\beta_{\mathrm{NW}} = \sum_{k_1 + k_2 = 2} \frac{1}{k_1! k_2!} \left[ \int s_1^{k_1} s_2^{k_2} K_2(s_1, s_2) \, ds_1 \, ds_2 \right]$$

$$\times \left\{ \frac{1}{f_2(t,z)} \frac{\partial^2}{\partial t^{k_1} \partial z^{k_2}} \alpha_1(t,z) - \frac{\mu(t,z)}{f_2(t,z)} \frac{\partial^2}{\partial t^{k_1} \partial z^{k_2}} f_2(t,z) \right\}$$

$$\times \tau_\mu \sqrt{\rho_\mu^{2k_2+1}},$$

$\Sigma_{\mathrm{NW}} = [\mathrm{Var}(Y|t,z)\|K_2\|^2]/f_2(t,z), \alpha_1(t,z) = \mu(t,z)f_2(t,z),$

and $f_2(t,z)$ is the joint density of $(T,Z)$.

THEOREM 3.2. *Under assumptions* A.3, A.5 *and* A.6, B.1–B.4, *and assuming* $\frac{h_{\mu,z}}{h_{\mu,t}} \to \rho_\mu$, *and* $nE(N)h_{\mu,t}^6 \to \tau_\mu^2$ *for some* $0 < \rho_\mu, \tau_\mu < \infty$,

$$\sqrt{n\bar{N}h_{\mu,t}h_{\mu,z}}[\hat{\mu}_L(t,z) - \mu(t,z)] \xrightarrow{\mathcal{D}} N(\beta_L, \Sigma_L),$$

*where*

$$\beta_L = \sum_{k_1+k_2=2} \frac{1}{k_1!k_2!} \left[ \int s_1^{k_1} s_2^{k_2} K_2(s_1,s_2)\, ds_1\, ds_2 \right] \frac{\partial^2}{\partial t^{k_1} \partial z^{k_2}} \mu(t,z)\tau_\mu \sqrt{\rho_\mu^{2k_2+1}},$$

$\Sigma_L = [\mathrm{Var}(Y|t,z)\|K_2\|^2]/f_2(t,z)$ *and* $f_2(t,z)$ *is the joint density of* $(T,Z)$.

*Asymptotic results for covariance functions.* We will need to consider three-dimensional smoothers to estimate the covariance function. Again, the asymptotic normalities of the Nadaraya–Watson kernel estimator and local linear estimator of the covariance function follow from Lemma D.2 in Appendix D. Here the Nadaraya–Watson kernel estimator of the covariance $\Gamma(t,s,z)$ is defined as

$$
\begin{aligned}
(3.3) \quad \hat{\Gamma}_{\mathrm{NW}}(t,s,z) = {} & \left( \sum_{i=1}^{n} \sum_{1 \le j \ne k \le N_i} K_3\left( \frac{t-T_{ij}}{h_{G,t}}, \frac{s-T_{ik}}{h_{G,t}}, \frac{z-Z_i}{h_{G,z}} \right) C_{ijk} \right) \\
& \times \left( \sum_{i=1}^{n} \sum_{1 \le j \ne k \le N_i} K_3\left( \frac{t-T_{ij}}{h_{G,t}}, \frac{s-T_{ik}}{h_{G,t}}, \frac{z-Z_i}{h_{G,z}} \right) \right)^{-1}.
\end{aligned}
$$

For notational convenience, we focus on the case of conventional kernel of order $(0,2)$ and denote

$$\sigma_i^2 = \iiint u_i^2 K_3(u_1,u_2,u_3)\, du_1\, du_2\, du_3$$

$$\text{for } i=1,2,3,$$

$$nE(N(N-1))h_{G,t}^6 h_{G,z} \to \tau_1^2, \qquad nE(N(N-1))h_{G,t}^2 h_{G,z}^5 \to \tau_2^2$$

and

$$v_3(t,s,z) = \mathrm{Var}((Y_1 - \mu(T_1,Z))(Y_2 - \mu(T_2,Z))|T_1=t, T_2=s, Z=z)$$

in the following two theorems.

THEOREM 3.3. *Under assumptions* A.4–A.6, B.5–B.8, *and assuming* $\frac{h_{G,z}}{h_{G,t}} \to \rho_G$ *and* $nE(N(N-1))h_{G,t}^7 \to \tau_G^2$ *for some* $0 < \rho_G, \tau_G < \infty$,

$$\sqrt{n\bar{N}(\bar{N}-1)h_{G,t}^2 h_{G,z}}\{\hat{\Gamma}_{\mathrm{NW}}(t,s,z) - \Gamma(t,s,z)\} \xrightarrow{\mathcal{D}} N(\gamma_{\mathrm{NW}}, \Omega_{\mathrm{NW}}),$$

*where*

$$\gamma_{\mathrm{NW}} = \frac{1}{2}\bigg\{ \sigma_1^2 \tau_1 \frac{d^2}{dt^2}\Gamma(t,s,z) + \sigma_2^2 \tau_1 \frac{d^2}{ds^2}\Gamma(t,s,z) + \sigma_3^2 \tau_2 \frac{d^2}{dz^2}\Gamma(t,s,z) \bigg\}$$
$$+ \bigg\{ \sigma_1^2 \tau_1 \bigg( \frac{d}{dt}\Gamma(t,s,z) \bigg) \bigg( \frac{d}{dt}g_3(t,s,z) \bigg)$$
$$+ \sigma_2^2 \tau_1 \bigg( \frac{d}{ds}\Gamma(t,s,z) \bigg) \bigg( \frac{d}{ds}g_3(t,s,z) \bigg)$$
$$+ \sigma_3^2 \tau_2 \bigg( \frac{d}{dz}\Gamma(t,s,z) \bigg) \bigg( \frac{d}{dz}g_3(t,s,z) \bigg) \bigg\} \bigg/ g_3(t,s,z),$$
$$\Omega_{\mathrm{NW}} = [v_3(t,s,z)\|K_3\|^2]/g_3(t,s,z)$$

*and* $g_3(t,s,z)$ *is the joint density of* $(T_1, T_2, Z)$.

THEOREM 3.4. *Under assumptions* A.4–A.6, B.5–B.8, *assuming* $\frac{h_{G,z}}{h_{G,t}} \to \rho_G$, *and* $nE(N(N-1))h_{G,t}^7 \to \tau_G^2$ *for some* $0 < \rho_G, \tau_G < \infty$,

$$\sqrt{n\bar{N}(\bar{N}-1)h_{G,t}^2 h_{G,z}}\{\hat{\Gamma}_L(t,s,z) - \Gamma(t,s,z)\} \xrightarrow{\mathcal{D}} N(\gamma_L, \Omega_L),$$

*where* $\gamma_L = \frac{1}{2}\{\sigma_1^2 \tau_1 \frac{d^2}{dt^2}\Gamma(t,s,z) + \sigma_2^2 \tau_1 \frac{d^2}{ds^2}\Gamma(t,s,z) + \sigma_3^2 \tau_2 \frac{d^2}{dz^2}\Gamma(t,s,z)\}$, $\Omega_L = [v_3(t,s,z)\|K_3\|^2]/g_3(t,s,z)$, *and* $g_3(t,s,z)$ *is the joint density of* $(T_1, T_2, Z)$.

*Remarks.* 1. The above asymptotic results reveal that standard optimal convergent rates for independent data are attained for all estimators when $E(N)$ is finite. For instance, the convergence rate for both the Nadaraya–Watson and local linear estimates for the mean function is $n^{1/3}$ which is the optimal convergence rate for a two-dimensional smoother under similar regularity conditions, and the convergence rate for both covariance function estimators is $n^{2/7}$, also optimal for a similar three-dimensional smoother.

2. The convergent rates of all estimators are faster when the expected number of measurements per subject $E(N) \to \infty$ as there are more data available per subject. For instance, the convergence rate for both mean function estimates and both covariance function estimates can be as arbitrarily close to $n^{2/5}$ when $E(N) \to \infty$. Note that $n^{2/5}$ is the optimal rate of convergence when the entire longitudinal process $Y(\cdot, z_i)$ can be observed for all subjects $i = 1, \ldots, n$; therefore smoothing is only required on the $z$-direction leading to a one-dimensional smoothing rate.

The asymptotic normality of the mFPCA covariance estimator can be handled similar to Theorem 3.2, but it is in fact much simpler if we follow the arguments of Theorem 2 in Yao (2007). The proof follows from the weak uniform convergence of $\hat{\mu}(t, z)$ in Lemma D.4, the asymptotic distributions of the estimators based on "raw covariances," $C_{ijk}$, are identical to those based on $\tilde{C}_{ijk} = \{Y_{ij} - \mu(T_{ij}, Z_i)\}\{Y_{ik} - \mu(T_{ik}, Z_i)\}$. Thus the Nadaraya–Watson estimator and local linear estimator of covariance based on $C_{ijk}$ are asymptotically equivalent to those estimators based on $\tilde{C}_{ijk}$. To save space, we present only the results for the local linear smoother in (2.8).

THEOREM 3.5. *Under assumptions* $h_{G^*} \to 0$, $nE(N^2)h_{G^*}^2 \to \infty$, $h_{G^*}^2 \times E(N^2) \to 0$, $nE(N(N-1))h_{G^*}^6 \to \tau^2$ *for some* $0 \le \tau < \infty$, A.5–A.6, *and* E.1–E.3,

$$\sqrt{n\bar{N}(\bar{N}-1)h_{G^*}^2}\{\hat{\Gamma}^*(t, s) - \Gamma^*(t, s)\} \xrightarrow{\mathcal{D}} N(\gamma^*, \Omega^*),$$

*where* $\gamma^* = \frac{\tau}{2} \int u^2 K_1(u)\, du \{\frac{d^2}{dt^2}\Gamma^*(t, s) + \frac{d^2}{ds^2}\Gamma^*(t, s)\}$, $\Omega^* = [v_2(t, s)\|K_1\|^4]/g_2(t, s)$, $v_2(t, s) = \mathrm{Var}((Y_1 - \mu(T_1, Z))(Y_2 - \mu(T_2, Z)))|T_1 = t, T_2 = s)$ *and* $g_2(t, s)$ *is the joint density of* $(T_1, T_2)$.

## 4. Simulation results.

We compare the performance of the two covariate adjusted FPCA approaches in Section 2 with the estimator in Yao, Müller and Wang (2005) which we term uFPCA with the prefix "u" suggesting that it is "unadjusted" FPCA. The simulation scheme is as follows: for each subject, a covariate $z$ is generated from $U(0, 1)$, its mean function is $\mu(t, z) = t + z \sin(t) + (1 - z) \cos(t)$ and its variance–covariance function is derived from two eigenfunctions $\phi_1(t, z) = -\cos(\pi(t + z/2))\sqrt{2}$ and $\phi_2(t, z) = \sin(\pi(t + z/2))\sqrt{2}$, for $0 \le t \le 1$ with eigenvalues $\lambda_1(z) = z/9$, $\lambda_2(z) = z/36$ and $\lambda_k = 0$ for $k \ge 3$. The specific principal component scores $A_{ik}(z)$ are generated from $N(0, \lambda_k(z))$, and the additional measurement errors are assumed to be normally distributed with mean 0 and variance $(0.05)^2$. For the measurement scheme $\{t_{ij}\}$ we use a nonequidistant "jittered" design. Specifically, an equally spaced grid $\{c_0, \ldots, c_{50}\}$ on $[0, 1]$ with $c_0 = 0$ and $c_{50} = 1$ is selected and jittered according to the plan $s_i = c_i + \epsilon_i$ where $\epsilon_i$ are i.i.d. with $N(0, 0.0001)$ and then constrained to be $s_i = 0$ if $s_i < 0$ and $s_i = 1$ if $s_i > 1$. Each curve is sampled at a random number of points, $\{t_{ij}\}$, $j = 1, \ldots, N_i$, where $N_i$ are chosen from a discrete uniform distribution $\{2, \ldots, 10\}$, and the locations of the measurements are randomly chosen from $\{s_1, \ldots, s_{49}\}$ without replacement. The simulation consists of 100 runs, and the number of subject is 100 in each run.

Epanechnikov kernels are used in the smoothing steps. The bandwidths for the mean surface estimator are chosen by leave-one-curve-out cross-validation while the bandwidths for the covariance estimator are chosen by a 10-fold cross-validation method to save computing time. Three criteria (AIC, BIC and fraction of variation explained (FVE) method) for choosing the value $K$ are also compared. The AIC and BIC are defined as in Yao, Müller and Wang (2005).

The FVE method is defined as the minimum number of components needed to explain at least a specified fraction of the total variation. In the simulation, we choose $K$ for uFPCA and mFPCA to be the minimum number $k$ satisfying $(\sum_{i=1}^{k} \lambda_i)/(\sum_{i=1} \lambda_i) \geq 0.80$, and for the fFPCA approach, this corresponds to selecting the smallest $k$ satisfying $\sum_{i=1}^{k} \lambda_i(z)/\sum_{i=1} \lambda_i(z) \geq 0.80$ for each subject with covariate value $z$. A major difference is that this type of FVE would allow subject-specific choice for the number of principal components in fFPCA. A problem is that the covariance estimate based on individually selected number of principal component may not yield a smooth covariance surface. To rectify this and to facilitate a uniform platform to compare the three approaches, we propose a global choice of $K$ based on the 90th percentile of the individually selected $k$'s for fFPCA. This global choice is somewhat objective and may give a slight benefit to fFPCA in fitting the observed data as compared to using either the mean or median value of $k$ as the global choice. Both AIC and BIC approaches tend to choose too many eigenfunctions so they can predict the data well, while FVE is the best for selecting the correct model. However, it is inferior to the others for prediction as evident in Table 2.

The mean integrated squared error of the covariance estimator for mFPCA is 0.00046, the biases and standard errors of the two eigenvalues are $-0.0102$ (s.d. $= 0.0121$) and $-0.0035$ (s.d. $= 0.0052$), respectively. The averaged estimated eigenfunction of the 100 simulations is close to the true eigenfunctions as shown in Figure 1. This suggests that the covariance estimator of mFPCA is sufficiently accurate. From Table 1 and Figure 2, the performance of fFPCA is generally satisfactory although the accuracy varies with the covariate. The estimate for the second eigenfunction at $Z = 0.1$ is poor due to the small eigenvalue 0.0028, so there is probably no need to include more than one eigenfunction for $Z = 0.1$.

Next, we compare the three different model selection criteria of choosing the number $K$ of eigenfunctions. We use the mean integrated squared error (MISE) for the true curves $X_i(t, z_i)$,

$$\text{(4.1)} \qquad \text{MISE} = \frac{1}{n} \sum_{i=1}^{n} \int_0^1 (X_i(t, z_i) - \hat{X}_i^K(t, z_i))^2 \, dt$$

as a criterion where $K$ is the number of eigenfunctions used to predict the trajectory of each subject. The corresponding mean squared fitting errors

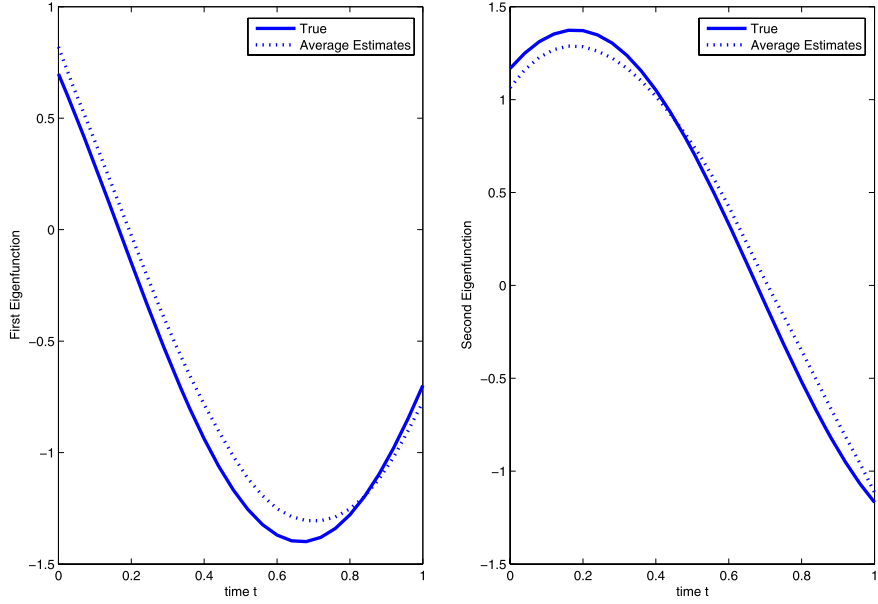

Fig. 1.    *First two eigenfunctions of the covariance and their estimates by mFPCA.*

(MSFE) is

$$\text{(4.2)} \qquad \text{MSFE} = \frac{1}{n} \sum_{i=1}^{n} \frac{1}{N_i} \sum_{j=1}^{N_i} (Y_{ij} - \hat{Y}_{ij})^2.$$

An outlier was detected in the 6th run for mFPCA predicted trajectory, so we include two results in Table 2, one with all simulations and one with this

TABLE 1

*Simulation results of fFPCA. The three rows corresponding to ISE are based on the average integrated squared errors of the 100 simulations where ISE of $\hat{g}(\cdot)$ for estimating a target function $g(\cdot)$ is defined as $\int_{\mathcal{T}} (\hat{g}(t) - g(t))^2 \, dt$. The rows corresponding to $\hat{\lambda}_i(z)$ are the biases and standard deviation (in bracket)*

| Covariate $z$ | 0.1 | 0.3 | 0.5 | 0.7 | 0.9 |
|---|---|---|---|---|---|
| ISE of $\hat{\Gamma}_L$ | 0.00015 | 0.00025 | 0.00071 | 0.0014 | 0.0030 |
| ISE of $\hat{\phi}_1(t, z)$ | 0.0294 | 0.0076 | 0.0071 | 0.0074 | 0.0112 |
| ISE of $\hat{\phi}_2(t, z)$ | 0.2720 | 0.0305 | 0.0242 | 0.0179 | 0.0300 |
| $\hat{\lambda}_1(z)$ | 0.0047 | −0.0041 | −0.0113 | −0.0202 | −0.0242 |
|  | (0.0073) | (0.0106) | (0.0181) | (0.0205) | (0.0333) |
| $\hat{\lambda}_2(z)$ | 0.0034 | 0.0001 | 0.0005 | −0.0002 | −0.0037 |
|  | (0.0045) | (0.0039) | (0.0057) | (0.0077) | (0.0094) |

TABLE 2
*Average MISE (4.1) and MSFE (4.2) in 100 simulation runs for the three approaches, the values in the parenthesis excludes one outlier occurred in the 6th run*

| | **MISE** | | | **MSFE** | | |
|---|---|---|---|---|---|---|
| | **FVE** | **AIC** | **BIC** | **FVE** | **AIC** | **BIC** |
| uFPCA | 0.0339 | 0.0215 | 0.0215 | 0.0047 | 0.0035 | 0.0036 |
| | (0.0325) | (0.0198) | (0.0197) | (0.0067) | (0.0065) | (0.0065) |
| mFPCA | 0.1075 | 0.0077 | 0.0076 | 0.0039 | 0.0024 | 0.0025 |
| | (0.0103) | (0.0063) | (0.0063) | (0.0050) | (0.0017) | (0.0017) |
| fFPCA | 0.0085 | 0.0077 | 0.0077 | 0.0039 | 0.0027 | 0.0027 |
| | (0.0085) | (0.0077) | (0.0077) | (0.0022) | (0.0015) | (0.0015) |

outlier run excluded. Not surprisingly, uFPCA is outperformed in general by both covariate adjusted approaches. When using the FVE method as a criterion of choosing $K$, fFPCA is slightly better than mFPCA. However, when using AIC or BIC as the criterion of choosing $K$, the performance of mFPCA is comparable to, if not better than, that of fFPCA. Consequently, if the purpose is to predict subject trajectories, mFPCA with BIC is recommended owing to its simplicity. For modeling purpose, fFPCA with the FVE method is preferred.

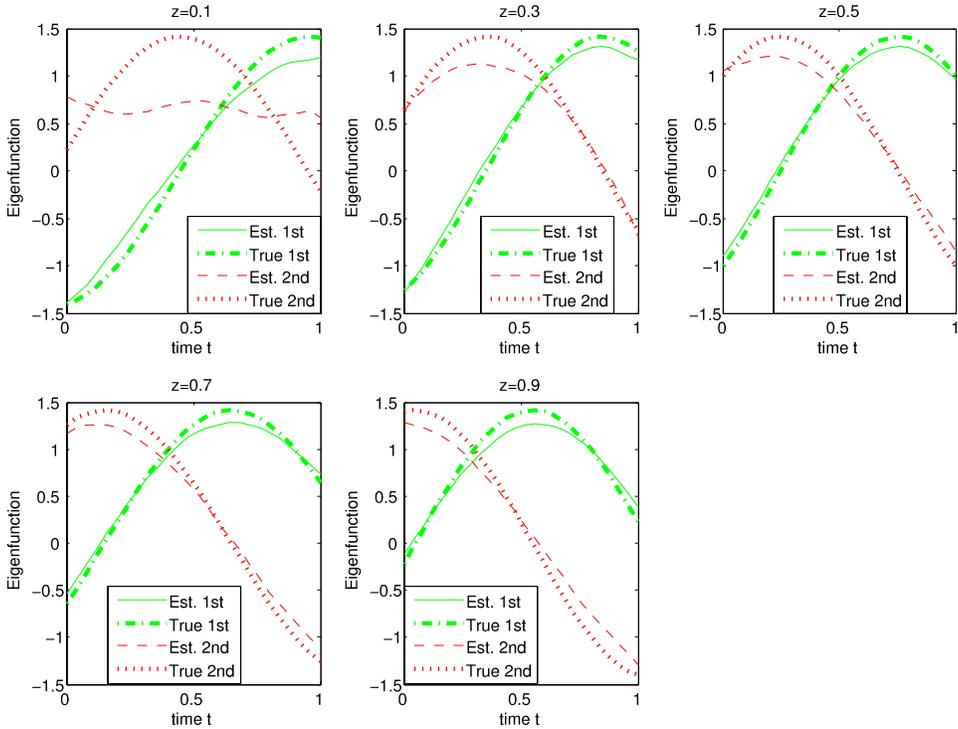

FIG. 2. *Means of the first two eigenfunctions estimated via fFPCA at five different values of the covarite.*

5. **Data application.** We illustrate the covariate-adjusted FPCA approaches through reproductive data for Mexican fruit flies. The study was conducted at the fruit fly mass rearing facility near Metapa, Chiapas, Mexico. Daily egg production (number of eggs) were recorded for a total of 1151 females [Carey et al. (2005)] till death. The goal here is to explore the influence of early reproduction, as measured by total reproduction by age 30 (in days), to reproduction pattern up to age 50. We exclude the infertile flies and those living less than 50 days. The latter is to provide a uniform platform to perform FPCA and consider only those who live at least around the average lifetime ($\approx$50.9 days) of fertile flies. Out of the remaining 567 flies, we further randomly selected 2 to 10 observations in the first 50 days, so we can compare the results for sparse data with the complete data to validate the new mFPCA and fFPCA approaches. In addition, we compare the new approaches with three different FPCA approaches that do not incorporate the covariate information. The first is the uFPCA in Yao, Müller and Wang (2005), the second is the reduced rank approach in James, Hastie and Suger (2000), termed rFPCA with "r" stands for reduced rank, and the third is a geometric approach in Peng and Paul (2009) similar to the reduced rank method but with a different algorithm. We term this approach "gFPCA" with "g" standing for geometric. Both rFPCA and gFPCA assume that $X(t)$ is a Gaussian processes, measurement errors are normally distributed, and use natural or B-spline bases to expand the eigenfunctions. Both approaches aim at maximizing the likelihood function, but rFPCA uses the EM algorithm to accomplish it and gFPCA tackles the likelihood functions directly with a Newton–Raphson method by exploiting the geometric structure of the eigenfunctions as they lie on a Stiefel manifold. As rFPCA serves as the initial estimates for gFPCA, the original code for rFPCA has been improved and included in an R package, fpca, which is available on the CRAN project.

As suggested in James, Hastie and Suger (2000), the number of bases in rFPCA is selected by 10-fold cross-validation likelihood and the number of eigenfunctions are reduced by the usual FVE (fraction of variation explained) method. For the Mexfly data, it selected 15 bases and the resulted numbers of eigenfunctions corresponding to 80% and 90% FVE, as reported in Table 3, are 9 and 11, respectively. The choice of the B-spline basis functions and the number of eigenfunctions for gFPCA are selected by a new cross-validated likelihood method proposed in Peng and Paul (2009) and they resulted in 8 bases and 5 eigenfunctions.

TABLE 3

*MSFEs of mFPCA, fFPCA, uFPCA and rFPCA with global K and the values in bracket are MSFEs based on sparse data*

| | **FVE (80%)** | | **FVE (90%)** | | **AIC** | | **BIC** | |
|---|---|---|---|---|---|---|---|---|
| | **MSFE** | $K$ | **MSFE** | $K$ | **MSFE** | $K$ | **MSFE** | $K$ |
| mFPCA | 614.1 (465.9) | 4 | 612.8 (447.9) | 6 | 611.8 (433.7) | 14 | 612.0 (436.4) | 10 |
| fFPCA | 614.9 (464.4) | 4 | 613.9 (454.4) | 5 | 612.8 (441.3) | 11 | 613.2 (445.7) | 7 |
| uFPCA | 684.6 (499.8) | 2 | 684.6 (499.8) | 2 | 680.8 (472.6) | 8 | 680.9 (473.6) | 6 |
| rFPCA | 720.2 (136.6) | 9 | 719.1 (131.5) | 11 | | | | |
| uFPCA | 681.0 (477.3) for $K = 4$, 680.8 (472.1) for $K = 10$, 680.7 (471.6) for $K = 14$ | | | | | | | |
| gFPCA | 785.1 (648.6) for $K = 5$ (based on the CV method), 784.8 (647.1) for $K = 6$ | | | | | | | |

Figure 3 shows the estimated mean surface of mFPCA and fFPCA for both sparse and complete data; it indicates that the daily reproductive rates are correlated with cumulative reproduction at young age, but that mean estimator works well even though data are sparse. Subjects with higher cumulative reproductions at a young age tend to have higher daily reproductive rates. Similar to the mean estimator, the covariance estimator of mFPCA also works very well when the data are sparse as Figure 4 shows. The estimated covariance function based on complete data is not so smooth as that based on sparse data because a smaller bandwidth was selected when there are substantially more data.

The mean square fitted errors for the five approaches are reported in Table 3. The performance of uFPCA, mFPCA and fFPCA are similar to those from the simulation study in Section 4; mFPCA is generally slightly better than fFPCA for sparse data, and both outperform uFPCA and gFPCA. The improvements of mFPCA and fFPCA over uFPCA appear marginal for sparse data, but this is due to the large measurement errors (the estimates of $\sigma$ by mFPCA, fFPCA, uFPCA are 25.34, 25.44, 24.81, respectively) present in the data. Since uFPCA only selects two eigenfunctions, we tried to check whether one can improve it by increasing the number of eigenfunctions. We use mFPCA as the gauge, and the lower portion of Table 3 reports additional results for uFPCA that utilizes the same number ($K = 4, 10$, and 14) of components as mFPCA. We also tried to include additional results for gFPCA to compare with mFPCA; however, the CV chose 8 bases and hence restricts $K$ to $K \leq 8$. This leads to only one additional case when $K = 6$ as the algorithm encountered singularity situation for the case with $K = 8$.

An intriguing phenomenon is the performance of rFPCA, which by far outperforms all other procedures for sparse data but not for the complete data where uFPCA, mFPCA and fFPCA all have smaller fitting errors.

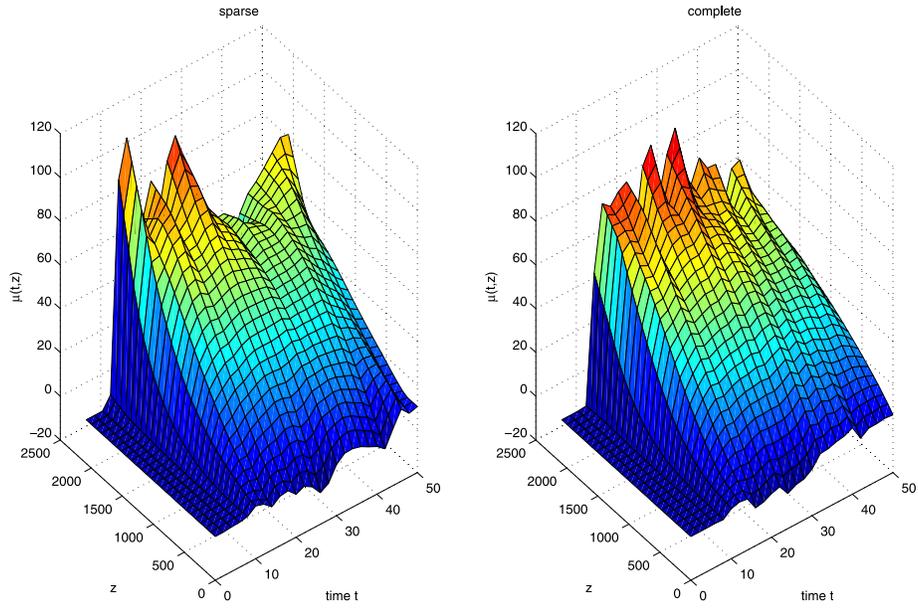

Fig. 3.  *Estimated mean surface for sparse and complete data.*

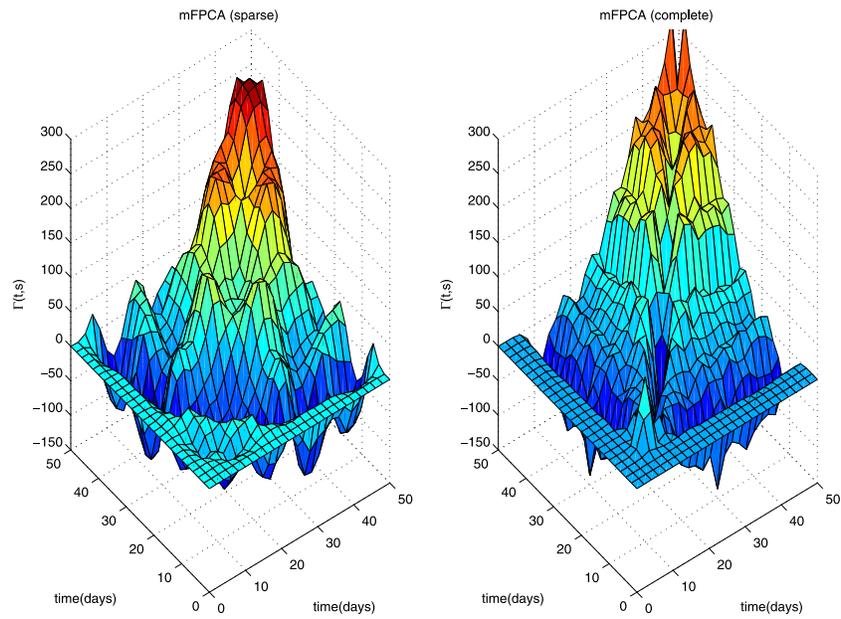

Fig. 4.  *Estimated covariance surfaces of mFPCA for sparse and compete data.*

This suggests an over-fitting problem and calls for further investigation. We tried to investigate this with simulations but could not reach any conclusion using the simulation setting in Section 4. Both algorithms in rFPCA and gFPCA encountered singularity situations or could not converge in many runs with the divergent problem more serious for gFPCA. It appears that the smoothing parameters for both methods are sensitive to the data.

In summary, this data supports the simpler covariate adjusted approach of adjusting just the mean but not the covariance. Additional benefit of the mFPCA approach is its computing speed. The computational time of the Mexican fruit flies data for fFPCA is 20 times more than mFPCA after the bandwidths for the mean and covariances functions have been selected. If we include the time to select those bandwidths, the gap is smaller as 10-fold CV was used to estimate the covariance functions for both mFPCA and fFPCA, leaving the leave-one-out CV for the mean function the most time-consuming part of the algorithm. However, the computational cost for fFPCA escalates as the total number of observations increases.

**6. Conclusions.** Through simulations and data analysis, we have shown that current approaches for functional principal component analysis may no longer be suitable for functional data when covariate information is available. Two alternatives are proposed to incorporate covariate effects on functional response data, adjusting the covariate effects on the mean function only (mFPCA) or adjusting the covariate effects for the covariance as well (fFPCA). Numerical evidence supports the simpler mean-adjusted approach especially when the purpose is to predict the trajectories $Y(t)$.

Besides the method itself, the criteria of choosing the number of eigenfunctions affect performance. Among the three criteria discussed in the simulation study, the FVE method based on the fraction of variation explained is more likely to pick the correct number of eigenfunctions than the other two criteria (AIC and BIC). When model fitting $X(t)$ is the main purpose of the data analysis, fFPCA with the FVE criterion is the best choice. However, fFPCA is time-consuming and mFPCA is just slightly less efficient than fFPCA in fitting $X(t)$ but could be more efficient than fFPCA in predicting $Y(t)$, so mFPCA might be an attractive approach to accommodate covariates.

Both FPCA approaches are model free and provide nonparametric estimates for both the fixed and random effects. The advantages of the principal-component based approach are: (1) less random effects are needed to fit the data; (2) it has the added value to reveal the modes of variation of the data and (3) it provides guidance to other parsimonious models such as a varying coefficient model or a linear mixed effects model. Developing formal inference procedures using either mFPCA or fFPCA approaches for model validation will be important future projects.

So far we have considered the theory and implementation for univariate covariate $Z$ only, but both can be extended to multivariate $Z$ conceptually and theoretically. The catch is the high-dimensional smoothing involved with a vector $Z$. Some dimension reduction on $Z$ will be needed for practical implementation and this will be another future research project.

## APPENDIX A: KERNEL FUNCTIONS

We consider both two and three-dimensional kernels that are symmetric with compact support. A kernel function $K_2: R^2 \to R$ is of order $(\boldsymbol{\nu}, \kappa)$ if

$$(A.1) \quad \iint u^{k_1} v^{k_2} K_2(u, v) \, du \, dv = \begin{cases} 0, & 0 \le k_1 + k_2 < \kappa, \\ & k_1 \ne \nu_1, k_2 \ne \nu_2, \\ (-1)^{|\boldsymbol{\nu}|} |\boldsymbol{\nu}|!, & k_1 = \nu_1, k_2 = \nu_2, \\ \ne 0, & k_1 + k_2 = \kappa, \end{cases}$$

where $\boldsymbol{\nu}$ is a multi-index $\boldsymbol{\nu} = (\nu_1, \nu_2)$ and $|\boldsymbol{\nu}| = \nu_1 + \nu_2$. A kernel function $K_3: R^3 \to R$ is of order $(\boldsymbol{\nu}, \kappa)$ if

$$
\begin{aligned}
(A.2) \quad & \iint u^{k_1} v^{k_2} w^{k_3} K_3(u, v, w) \, du \, dv \, dw \\
& = \begin{cases} 0, & 0 \le \sum_{i=1}^{3} k_i < \kappa, k_i \ne \nu_i, \\ & \text{for } i = 1, 2, 3, \\ (-1)^{|\boldsymbol{\nu}|} |\boldsymbol{\nu}|!, & k_1 = \nu_1, k_2 = \nu_2, k_3 = \nu_3, \\ \ne 0, & k_1 + k_2 + k_3 = \kappa, \end{cases}
\end{aligned}
$$

where $\boldsymbol{\nu}$ is a multi-index $\boldsymbol{\nu} = (\nu_1, \nu_2, \nu_3)$ and $|\boldsymbol{\nu}| = \nu_1 + \nu_2 + \nu_3$.

## APPENDIX B: ASSUMPTIONS

For bandwidth sequences $h_1 = h_1(n)$ and $h_2 = h_2(n)$, the notation $h_1 \asymp h_2$ means they are of the same order, and, namely, $h_1/h_2$ stays away from 0 and $\infty$. We denote $h_{\mu,t}$ and $h_{\mu,z}$ as the two bandwidth sequences for the mean function estimator in the coordinates $\mathcal{T}$ and $\mathcal{Z}$. Similarly, $h_{G,t}$ and $h_{G,z}$ are the two bandwidth sequences for the covariance estimator. The assumptions on the bandwidths are listed in A.1–A.4; the assumptions of the measurement schedule are in A.5, A.6 and A.7 is a common assumption while the covariance is estimated. The bandwidth assumptions and the measurement schedule assumptions are required to show that the local property of the kernel-based estimators holds for longitudinal or functional data with the presence of within-subject correlation.

A.1 $h_{\mu,t} \asymp h_{\mu,z} \asymp h$, $h \to 0$, $nh^{|\nu|+2} \to \infty$, and $nh^{2\kappa+2} < \infty$.

A.2 $h_{G,t} \asymp h_{G,z} \asymp h$, $h \to 0$, $nh^{|\nu|+3} \to \infty$, and $nh^{2\kappa+3} < \infty$.

A.3 $h_{\mu,t} \asymp h_{\mu,z} \asymp h$, $h \to 0$, $nE(N)h^{|\nu|+2} \to \infty$, $E(N)h \to 0$ and $nE(N) \times h^{2\boldsymbol{\kappa}+2} < \infty$.

A.4 $h_{G,t} \asymp h_{G,z} \asymp h$, $h \to 0$, $nE(N^2)h^{|\nu|+3} \to \infty$, $E(N^2)h^2 \to 0$ and $nE(N \times (N-1))h^{2\boldsymbol{\kappa}+3} < \infty$.

A.5 The number of observations $N_i(n)$ for the $i$th subject is a random variable with $N_i(n) \sim N(n)$ where $N(n)$ is a positive integer-valued random variable with $\limsup_{n\to\infty} \frac{EN(n)^2}{[EN(n)]^2}$ and $\limsup_{n\to\infty} \frac{EN(n)^4}{(EN(n)^2)^2}$ both finite. Moreover, $N_i(n), i = 1, \ldots, n$ are i.i.d.

A.6 The observational times $T_{ij}$ and measurements $Y_{ij}$ are independent of the number of measurements $N(n)$.

A.7 $E\{(Y - \mu(T, Z))^4\} < \infty$.

For random design, we assume $(T_{ij}, Z_i, Y_{ij})$ have the same distribution as $(T, Z, Y)$ with joint p.d.f. $f_3(t, z, y)$, and the observation times $T_{ij}$ are i.i.d. with p.d.f. $f(t)$, but dependency is allowed among observations from the same subject. The joint p.d.f.'s of $(T, Z)$, $(T_1, T_2, Y_1, Y_2)$, $(T_1, T_2, Z, Y_1, Y_2)$, $(T_1, T_2, T_1', T_2', Y_1, Y_2, Y_1', Y_2')$, $(T_1, T_2, T_1', T_2', Z, Y_1, Y_2, Y_1', Y_2')$, $(T_1, T_2)$, and $(T_1, T_2, Z)$ are, respectively, $f_2(t, z)$, $f_4(t_1, t_2, y_1, y_2)$, $f_5(t_1, t_2, z, y_1, y_2)$, $f_8(t_1, t_2, t_1', t_2', y_1, y_2, y_1', y_2')$, $f_9(t_1, t_2, t_1', t_2', z, y_1, y_2, y_1', y_2')$, $g_2(t_1, t_2)$ and $g_3(t_1, t_2, z)$. The following regularity conditions for marginal or joint p.d.f.'s and mean or covariance functions along with the bandwidths assumptions are used to show the consistency results in Section 3 and Appendices C and D. B.1–B.4 are for the asymptotic results of two-dimensional kernel-based estimators while B.5–B.8 are for the asymptotic results of three-dimensional kernel-based estimators.

The following type of continuity, as defined in Yao (2007), will be needed:

DEFINITION 1. A real function $f(x, y) : R^{n+m} \to R$ is continuous on $A \subseteq R^n$ uniformly in $y \in R^m$, if given any $x \in A$ and $\varepsilon > 0$ there exists a neighborhood of $x$ not depending on $y$, say $U(x)$, such that $|f(x', y) - f(x, y)| < \varepsilon$ for all $x' \in U(x)$ and $y \in R^m$.

B.1 $\frac{d^{\boldsymbol{\kappa}}}{dt^{k_1} dz^{k_2}} f_2(t, z)$ exists and is continuous on $\{(t, z)\}$ for $k_1 + k_2 = \boldsymbol{\kappa}$, $0 \le k_1, k_2 \le \boldsymbol{\kappa}$, and $f_2(t, z) > 0$.

B.2 $f_3(t, z, y)$ is continuous on $\{(t, z)\}$ uniformly in $y \in R$; $\frac{d^{\boldsymbol{\kappa}}}{dt^{k_1} dz^{k_2}} f_3(t, z, y)$ exists and is continuous on $\{(t, z)\}$ uniformly in $y \in R$, for $k_1 + k_2 = \boldsymbol{\kappa}$, $0 \le k_1, k_2 \le \boldsymbol{\kappa}$.

B.3 $f_5(t_1, t_2, z, y_1, y_2)$ is continuous on $\{(t_1, t_2, z)\}$ uniformly in $(y_1, y_2) \in R^2$.

B.4 $\frac{d^{\boldsymbol{\kappa}}}{dt^{k_1} dz^{k_2}} \mu(t, z)$ exists and is continuous on $\{(t, z)\}$, for $k_1 + k_2 = \boldsymbol{\kappa}$, $0 \le k_1, k_2 \le \boldsymbol{\kappa}$.

B.5 $\frac{d^{\boldsymbol{\kappa}}}{dt^{k_1} ds^{k_2} dz^{k_3}} g_3(t, s, z)$ exists and is continuous on $\{(t, s, z)\}$ for $k_1 + k_2 + k_3 = \boldsymbol{\kappa}$, $0 \le k_1, k_2, k_3 \le \boldsymbol{\kappa}$, and $g_3(t, s, z) > 0$.

B.6 $f_5(t, s, z, y_1, y_2)$ is continuous on $\{(t, s, z)\}$ uniformly in $(y_1, y_2) \in R^2$; $\frac{d^{\boldsymbol{\kappa}}}{dt^{k_1} ds^{k_2} dz^{k_3}} f_5(t, s, z, y_1, y_2)$ exists and is continuous on $\{(t, s, z)\}$ uniformly in $(y_1, y_2)$, for $k_1 + k_2 + k_3 = \boldsymbol{\kappa}$, $0 \leq k_1, k_2, k_3 \leq \boldsymbol{\kappa}$.

B.7 $f_9(t_1, t_2, t_1', t_2', z, y_1, y_2, y_1', y_2')$ is continuous on $\{(t_1, t_2, t_1', t_2', z)\}$ uniformly in $(y_1, y_2, y_1', y_2') \in R^4$.

B.8 $\frac{d^{\boldsymbol{\kappa}}}{dt^{k_1} ds^{k_2} dz^{k_3}} \Gamma(t, s, z)$ exists and is continuous on $\{(t, s, z)\}$, for $k_1 + k_2 + k_3 = \boldsymbol{\kappa}$, $0 \leq k_1, k_2, k_3 \leq \boldsymbol{\kappa}$.

## APPENDIX C: PROOFS OF THEOREMS 3.1 AND 3.2

Given an integer $Q \geq 1$ and for $q = 1, \ldots, Q$, let $\psi_q : R^3 \to R$ satisfy:

C.1 $\psi_q(t, z, y)$'s are continuous on $U(\{t, z\})$ uniformly in $y \in R$;

C.2 The functions $\frac{\partial^p}{\partial t^{p_1} \partial z^{p_2}} \psi_q(t, z, y)$ exist for all arguments $(t, z, y)$ and are continuous on $U(\{t, z\})$ uniformly in $y \in R$, for $p_1 + p_2 = p$ and $0 \leq p_1, p_2 \leq p$.

The kernel-weighted averages for two-dimensional smoothers are defined as

(C.1)    $\Psi_{qn} = \dfrac{1}{nEN h_{\mu,t}^{\nu_1+1} h_{\mu,z}^{\nu_2+1}} \sum_{i=1}^{n} \sum_{j=1}^{N_i} \psi_q(T_{ij}, Z_i, Y_{ij}) K_2\left(\dfrac{t - T_{ij}}{h_{\mu,t}}, \dfrac{z - Z_i}{h_{\mu,z}}\right),$

where $K_2$ is a kernel function of order $(\nu, \kappa)$ [defined in (A.1)], $h_{\mu,t}$, and $h_{\mu,z}$ are bandwidths associated with $t$ and $z$, respectively. We will see later that the Nadaraya–Watson estimator and local linear estimator each involves two and four such $\psi_q$ functions yielding $Q = 2$ and 4, respectively. Let

$$\alpha_q(t, z) = \frac{\partial^{|\boldsymbol{\nu}|}}{\partial t^{\nu_1} \partial z^{\nu_2}} \int \psi_q(t, z, y) f_3(t, z, y) \, dy$$

and

$$\sigma_{qr}(t, z) = \int \psi_q(t, z, y) \psi_r(t, z, y) f_3(t, z, y) \, dy \, \|K_2\|^2,$$

where $f_3(t, z, y)$ is the joint density of $(T, Z, Y)$, $\|K_2\|^2 = \int K_2^2$ and $1 \leq q, r \leq Q$.

We first provide the asymptotic normality of kernel-weighted averages for two-dimensional smoothers based on longitudinal data. Lemma C.1 extends Theorem 4.1 of Bhattacharya and Müller (1993) from a univariate smoother on independent observations to a bivariate smoother on correlated longitudinal observations. Lemma C.2 provides the key steps for the asymptotic results of the Nadaray–Watson and local linear estimators.

LEMMA C.1. *Under assumptions* A.3, A.5 *and* A.6, B.1–B.4, *and* C.1 *and* C2,

$$
\text{(C.2)} \quad \begin{aligned} &\sqrt{nENh_{\mu,t}^{2\nu_1+1}h_{\mu,z}^{2\nu_2+1}}[(\Psi_{1n},\ldots,\Psi_{Qn})^T - (E\Psi_{1n},\ldots,E\Psi_{Qn})^T] \\ &\qquad \xrightarrow{\mathcal{D}} N(0,\Sigma). \end{aligned}
$$

PROOF. We will show this through Cramér–Wold device and Lindeberg CLT. Let

$$
A = \sqrt{nENh_{\mu,t}^{2\nu_1+1}h_{\mu,z}^{2\nu_2+1}}\sum_{q=1}^{Q}a_q[\Psi_{qn} - E(\Psi_{qn})],
$$

where $a_q$, $1 \le q \le Q$, are given constants. Observing that $A = \sum_{i=1}^{n}U_i$ where

$$
U_i = \frac{1}{\sqrt{nENh_{\mu,t}h_{\mu,z}}}\sum_{q=1}^{Q}\sum_{j=1}^{N_i}a_q\psi_q(T_{ij},Z_i,Y_{ij})K_2\left(\frac{t-T_{ij}}{h_{\mu,t}},\frac{z-Z_i}{h_{\mu,z}}\right)
$$

$$
- \sum_{q=1}^{Q}\frac{a_q}{n}\sqrt{nENh_{\mu,t}^{2\nu_1+1}h_{\mu,z}^{2\nu_2+1}}E\Psi_{qn}
$$

and $U_i$'s are i.i.d. mean zero random variables. To verify the Lindeberg condition, we need $\text{Var}(U_i)$, $1 \le i \le n$. First, we show

$$
nENh_{\mu,t}^{2\nu_1+1}h_{\mu,z}^{2\nu_2+1}\text{cov}(\Psi_{qn},\Psi_{rn}) = \sigma_{qr}(t,z) + o(1).
$$

To see this, consider $nENh_{\mu,t}^{2\nu_1+1}h_{\mu,z}^{2\nu_2+1}\text{cov}(\Psi_{qn},\Psi_{rn}) = I_1 - I_2$ where

$$
\begin{aligned} I_1 = \frac{1}{h_{\mu,t}h_{\mu,z}}E\Bigg[&\frac{1}{EN}\left\{\sum_{j=1}^{N}\psi_q(T_j,Z,Y_j)K_2\left(\frac{t-T_j}{h_{\mu,t}},\frac{z-Z}{h_{\mu,z}}\right)\right\} \\ &\times\left\{\sum_{l=1}^{N}\psi_q(T_l,Z,Y_l)K_2\left(\frac{t-T_l}{h_{\mu,t}},\frac{z-Z}{h_{\mu,z}}\right)\right\}\Bigg] \end{aligned}
$$

and

$$
\begin{aligned} I_2 = \frac{EN}{h_{\mu,t}h_{\mu,z}}E\Bigg[&\frac{1}{EN}\sum_{j=1}^{N}\psi_q(T_j,Z,Y_j)K_2\left(\frac{t-T_j}{h_{\mu,t}},\frac{z-Z}{h_{\mu,z}}\right)\Bigg] \\ &\times E\Bigg[\frac{1}{EN}\sum_{l=1}^{N}\psi_r(T_l,Z,Y_l)K_2\left(\frac{t-T_l}{h_{\mu,t}},\frac{z-Z}{h_{\mu,z}}\right)\Bigg]. \end{aligned}
$$

It is obvious that $I_2 = o(1)$. As for $I_1$, it can be decomposed to $I_1 = Q_1 + Q_2$ where

$$Q_1 = \frac{1}{h_{\mu,t} h_{\mu,z}} E\left[\frac{1}{EN}\left(\sum_{j=1}^{N} \psi_q(T_j, Z, Y_j)\psi_r(T_j, Z, Y_j)K_2^2\left(\frac{t-T_j}{h_{\mu,t}}, \frac{z-Z}{h_{\mu,z}}\right)\right)\right]$$

$$= \frac{1}{h_{\mu,t} h_{\mu,z}} E\left[\psi_q(T, Z, Y)\psi_r(T, Z, Y)K_2^2\left(\frac{t-T}{h_{\mu,t}}, \frac{z-Z}{h_{\mu,z}}\right)\right]$$

$$= \sigma_{qr}^2(t, z) + o(1)$$

and

$$Q_2 = \frac{h_{\mu,t} E(N(N-1))}{EN} E\{\psi_q(t - t_1 h_{\mu,t}, z - h_{\mu,z} s_2, y_1)$$
$$\times \psi_r(t - t_2 h_{\mu,t}, z - h_{\mu,z} s_2, y_2)K_2(t_1, s_2)K_2(t_2, s_2)\}$$

$$= (h_{\mu,t} EN)\frac{E(N(N-1))}{(EN)^2}$$

$$\times \left\{\int \psi_q(t, z, y_1)\psi_r(t, z, y_2)f_5(t, t, z, y_1, y_2)\, dy_1\, dy_2\right.$$

$$\left.\times \left(\int K_2(t_1, s_2)K_2(t_2, s_2)\, dt_1\, dt_2\, ds_2\right) + o(h)\right\} = o(1).$$

Therefore, we can have $\mathrm{Var}(U_i) = \frac{1}{n}(a^T \Sigma a + o(1))$ where $a^T = (a_1, \ldots, a_Q)$. Let $B_n = \sum_{i=1}^{n} \mathrm{Var}(U_i) = a^T \Sigma a + o(1)$. In order to apply Lindeberg CLT, we need to prove

$$\lim_{n\to\infty} \frac{1}{B_n^2}\sum_{i=1}^{n} E[U_i^2 1_{\{|U_i|>\epsilon B_n\}}] = 0 \qquad \forall \epsilon > 0,$$

where $1_{\{\cdot\}}$ is an indicator function and it suffices to prove

$$\lim_{n\to\infty} nE[U_1^2 1_{\{U_1^2 > \epsilon^2 B_n^2\}}] = 0.$$

Using the fact $(a+b)^2 \le 2a^2 + 2b^2$, we can get

$$nE\{U_1^2 1_{\{U_1^2 > \epsilon^2 B_n^2\}}\}$$

$$\le 2nE\left\{\frac{1}{nENh_{\mu,t}h_{\mu,z}}\left[\sum_{q=1}^{Q}\sum_{j=1}^{N_1} a_q\psi_q(T_{1j}, Z_1, Y_{1j})\right.\right.$$

$$\left.\left.\times K_2\left(\frac{t-T_{1j}}{h_{\mu,t}}, \frac{z-Z_1}{h_{\mu,z}}\right)\right]^2 1_{\{\eta\}}\right\} + o(1),$$

where

$$\eta = \frac{1}{ENh_{\mu,t}h_{\mu,z}} \left( \sum_{q=1}^{Q} \sum_{j=1}^{N_1} a_q \psi_q(T_{1j}, Z_1, Y_{1j}) K_2 \left( \frac{t - T_{1j}}{h_{\mu,t}}, \frac{z - Z_1}{h_{\mu,z}} \right) \right)^2$$

$$> \frac{n\epsilon^2}{2} (a^T \Sigma a + o(1)) - o(1).$$

Observing that the term $[\sum_{q=1}^{Q} \sum_{j=1}^{N_1} a_q \psi_q(T_{1j}, Z_1, Y_{1j}) K_2(\frac{t-T_{1j}}{h_{\mu,t}}, \frac{z-Z_1}{h_{\mu,z}})]^2$ is dominated by $\sum_{q=1}^{Q} \sum_{r=1}^{Q} \sum_{j=1}^{N_1} a_q a_r \psi_q(T_{1j}, Z_1, Y_{1j}) \psi_r(T_{1j}, Z_1, Y_{1j}) \times K_2^2(\frac{t-T_{1j}}{h_{\mu,t}}, \frac{z-Z_1}{h_{\mu,z}})$, and using change of variable, we arrive at

$$nE\{U_1^2 1_{\{U_1^2 > \epsilon^2 B_n^2\}}\} \le 2E\{\Upsilon 1_{\{\Upsilon > n\epsilon^2/2(a^T \Sigma a + o(1)) - o(1)\}}\} + o(1),$$

where $\Upsilon = \sum_{q=1}^{Q} \sum_{r=1}^{Q} a_q a_r \psi_q(t - s_1 h_{\mu,t}, z - s_2 h_{\mu,z}, Y) \psi_r(t - s_1 h_{\mu,t}, z - s_2 h_{\mu,z}, Y) K_2^2(s_1, s_2)$, $\frac{t-T}{h_{\mu,t}} = s_1$, and $\frac{z-Z}{h_{\mu,z}} = s_2$. So far, we have shown that $\lim_{n\to\infty} \frac{n\epsilon^2}{2}(a^T \Sigma a + o(1)) = \infty$ for any given $\epsilon > 0$. This implies that Lindeberg condition holds and the proof of the lemma is thus complete. $\quad\square$

LEMMA C.2. *Let $H: R^Q \to R$ be a function with continuous first-order derivatives, $DH(v) = (\frac{\partial}{\partial x_1} H(v), \dots, \frac{\partial}{\partial x_Q} H(v))^T$, and $\bar{N} = \frac{1}{n} \sum_{i=1}^{n} N_i$. Under assumptions* A.3, A.5 *and* A.6, B.1–B.4, C.1 *and* C.2, *and assuming $\frac{h_{\mu,z}}{h_{\mu,t}} \to \rho_\mu$ and $nE(N)h_{\mu,t}^{2\kappa+2} \to \tau_\mu^2$ for some $0 < \rho_\mu, \tau_\mu < \infty$,*

$$\sqrt{n\bar{N}h_{\mu,t}^{2\nu_1+1}h_{\mu,z}^{2\nu_2+1}}[H(\Psi_{1n}, \dots, \Psi_{Qn}) - H(\alpha_1, \dots, \alpha_Q)]$$

$$\xrightarrow{\mathcal{D}} N(\beta_H, [DH(\alpha_1, \dots, \alpha_Q)]^T \Sigma [DH(\alpha_1, \dots, \alpha_Q)]),$$

*where $\Sigma = (\sigma_{qr})_{1 \le q, r \le l}$, and*

$$\beta_H = \sum_{k_1 + k_2 = \kappa} \frac{(-1)^{\kappa}}{k_1! k_2!} \left[ \int s_1^{k_1} s_2^{k_2} K_2(s_1, s_2) \, ds_1 \, ds_2 \right]$$

$$\times \left\{ \sum_{q=1}^{Q} \frac{\partial H}{\partial \alpha_q} [(\alpha_1, \dots, \alpha_Q)^T] \frac{\partial^{k_1 + k_2 - \nu_1 - \nu_2}}{\partial t^{k_1 - \alpha_q} \partial z^{k_2 - \nu_2}} \alpha_q(t, z) \right\} \tau_\mu \sqrt{\rho_\mu^{2k_2+1}}.$$

PROOF. It suffices to show this theorem with $\bar{N}$ replaced by $E(N)$ due to *Slutsky theorem*. We first handle the asymptotic bias term by showing that

$$(C.3) \qquad \sqrt{nE(N)h_{\mu,t}^{2\nu_1+1}h_{\mu,z}^{2\nu_2+1}}[H(\Psi_{1n}, \dots, \Psi_{Qn}) - H(\alpha_1, \dots, \alpha_Q)] \to \beta_H.$$

By conditioning on the value of $N$, it is easy to see that

$$E(\Psi_{qn}) = \frac{1}{h_{\mu,t}^{\nu_1+1} h_{\mu,z}^{\nu_2+1}} E\left(\psi_q(T, Z, Y) K_2\left(\frac{t-T}{h_{\mu,t}}, \frac{z-Z}{h_{\mu,z}}\right)\right).$$

Let $\frac{t-T}{h_{\mu,t}} = s_1$ and $\frac{z-Z}{h_{\mu,z}} = s_2$; it follows from Taylor's expansion of order $|k|$ on $\Psi_{qn}$'s and Taylor's expansion on $H$ that

$$E(\Psi_{qn}) = \alpha_q(t, z) + \sum_{k_1+k_2=|k|} \frac{(-1)^{|k|}}{k_1! k_2!} \left[\iint s_1^{k_1} s_2^{k_2} K_2(s_1, s_2)\, ds_1\, ds_2\right]$$

$$\times \left[\frac{\partial^{|k|-(\nu_1+\nu_2)}}{\partial t^{k_1-\nu_1} \partial z^{k_2-\nu_2}} \alpha_q(t, z)\right] h_{\mu,t}^{k_1-\nu_1} h_{\mu,z}^{k_2-\nu_2}$$

$$+ o(h_{\mu,t}^{k_1-\nu_1}) + o(h_{\mu,z}^{k_2-\nu_2}).$$

Combining Lemma C.1, the continuity of $DH$ at $(\alpha_1, \ldots, \alpha_Q)^T$ and the $\delta$-method, we have

$$\begin{aligned}
\text{(C.4)} \quad & \sqrt{nENh_{\mu,t}^{2\nu_1+1} h_{\mu,z}^{2\nu_2+1}} [H(\Psi_{1n}, \ldots, \Psi_{Qn}) - H(E\Psi_{1n}, \ldots, E\Psi_{Qn})] \\
& \xrightarrow{\mathcal{D}} N(0, [DH(\alpha_1, \ldots, \alpha_Q)]^T \Sigma [DH(\alpha_1, \ldots, \alpha_Q)]).
\end{aligned}$$

The lemma now follows from (C.3) and (C.4).  □

PROOF OF THEOREM 3.1.  Let $\psi_1(u_1, u_2, u_3) = u_3$, $\psi_2(u_1, u_2, u_3) = 1$, and $H(x_1, x_2) = x_1/x_2$, then $\hat{\mu}_{\text{NW}} = H(\Psi_{1n}, \Psi_{2n})$, $DH(\alpha_1, \alpha_2) = (1/\alpha_2, -\alpha_1/\alpha_2^2)$, $\alpha_1(t, z) = \mu(t, z) f_2(t, z)$, $\alpha_2(t, z) = f_2(t, z)$. Applying the results of Lemma C.2, the bias is

$$\begin{aligned}
\beta_{\text{NW}} = \sum_{k_1+k_2=2} \frac{1}{k_1! k_2!} & \left[\int s_1^{k_1} s_2^{k_2} K_2(s_1, s_2)\, ds_1\, ds_2\right] \\
& \times \left\{\frac{1}{\alpha_2(t, z)} \frac{\partial^2}{\partial t^{k_1} \partial z^{k_2}} \alpha_1(t, z) - \frac{\alpha_1(t, z)}{(\alpha_2(t, z))^2} \frac{\partial^2}{\partial t^{k_1} \partial z^{k_2}} \alpha_2(t, z)\right\} \\
& \times \tau_\mu \sqrt{\rho_\mu^{2k_2+1}}
\end{aligned}$$

and the components of $\Sigma$ are

$$\sigma_{11}(t, z) = [\text{Var}(Y | T = t, Z = z) + \mu^2(t, z)] f_2(t, z) \|K_2\|^2,$$

$$\sigma_{12}(t, z) = \sigma_{21}^2(t, z) = \mu(t, z) f_2(t, z) \|K_2\|^2, \qquad \sigma_{22}(t, z) = f_2(t, z) \|K_2\|^2.$$

Therefore, $\Sigma_{\text{NW}} = \frac{\text{Var}(Y | t, z)}{f_2(t, z)} \|K_2\|^2$, and the result follows.  □

The next lemma follows similar arguments as in Lemma 1 of Yao, Müller and Wang (2005) except that a two-dimensional Fourier transformation is employed here whereas their arguments involve only a one-dimensional Fourier transformation. Define the Fourier transforms of $K_2(u,v)$ and $K_3(u_1, u_2, u_3)$ by

$$\zeta_1(t,z) = \iint \exp(-(iut + iwz))K_2(u,w)\,du\,dw$$

and

$$\zeta_2(t,s,z) = \iiint \exp(-(iu_1 t + iu_2 s + iu_3 z))K_3(u_1, u_2, u_3)\,du_1\,du_2\,du_3.$$

They satisfy:

D.1 $\zeta_1(t,z)$ is absolutely integrable;
D.2 $\zeta_2(t,s,z)$ is absolutely integrable.

LEMMA C.3. *Under assumptions* A.1, A.5–A.7, B.1–B.4, C.1 *and* C.2 *and* D.1,

$$\sup_{t \in \mathcal{T}; z \in \mathcal{Z}} |\Psi_{qn} - \alpha_q| = O_p\left(\frac{1}{\sqrt{n}h_\mu^{|\nu|+2}}\right) \qquad \text{where } h_\mu \asymp h_{\mu,t} \asymp h_{\mu,z}.$$

PROOF. Since

$$E\left\{\sup_{t \in \mathcal{T}; z \in \mathcal{Z}} |\Psi_{qn} - \alpha_q|\right\} \leq E\left\{\sup_{t \in \mathcal{T}; z \in \mathcal{Z}} |\Psi_{qn} - E(\Psi_{qn})|\right\}$$
$$+ \sup_{t \in \mathcal{T}; z \in \mathcal{Z}} |E(\Psi_{qn}) - \alpha_q|$$

and Taylor's expansion implies $E(\Psi_{qn}) = \alpha_q + O(h_\mu^{k-\nu_1-\nu_2}) = \alpha_q + O(\frac{1}{\sqrt{n}h_\mu^{|\nu|+2}})$. It remains to show the correct order of the first term. To this end, we employ the inverse Fourier transformation

$$\zeta_1(t,z) = \iint \exp(-iut - iwz)K_2(u,w)\,du\,dw,$$

which implies

$$K_2\left(\frac{t - T_{\ell j}}{h_{\mu,t}}, \frac{z - Z_\ell}{h_{\mu,z}}\right)$$
$$= \left(\frac{1}{2\pi}\right)^2 \iint \exp\left(iu\left(\frac{t - T_{\ell j}}{h_{\mu,t}}\right) + iw\left(\frac{z - Z_\ell}{h_{\mu,z}}\right)\right)\zeta_1(u,w)\,du\,dw.$$

Let $\varphi_{qn}(u,w) = \frac{1}{n}\sum_{\ell=1}^{n}\frac{1}{E(N)}\sum_{j=1}^{N_{\ell}}\exp(iuT_{\ell j}+iwZ_{\ell})\psi(T_{\ell j},Z_{\ell},Y_{\ell j})$, and by plugging equation (C.5) into $\Psi_{qn}$, we obtain

$$\Psi_{qn} = \left(\frac{1}{2\pi}\right)^2 \frac{1}{h_{\mu,t}^{\nu_1}h_{\mu,z}^{\nu_2}}\int\int \varphi_{qn}(u,w)\exp(-iut-iwz)\zeta_1(h_{\mu,t}u,h_{\mu,z}w)\,du\,dw.$$

Therefore,

$$\sup_{t\in\mathcal{T};z\in\mathcal{Z}}|\Psi_{qn}-E(\Psi_{qn})|$$

$$\leq \left(\frac{1}{2\pi}\right)^2 \frac{1}{h_{\mu,t}^{\nu_1}h_{\mu,z}^{\nu_2}}$$

$$\times\int\int|\varphi_{qn}(u,w)-E(\varphi_{qn}(u,w))||\zeta_1(h_{\mu,t}u,h_{\mu,z}w)|\,du\,dw.$$

By the facts that $E(|\varphi_{qn}(u,w)-E(\varphi_{qn}(u,w))|) \leq (E\{\varphi_{qn}(u,w)-E(\varphi_{qn}(u,w))\}^2)^{1/2}$ and $\{\mathbf{T}_{\boldsymbol{\ell}},\mathbf{Y}_{\boldsymbol{\ell}},N_{\ell}\}$ are i.i.d. where $\mathbf{T}_{\boldsymbol{\ell}}=(T_{\ell 1},\ldots,T_{\ell N_{\ell}})^T$ and $\mathbf{Y}_{\boldsymbol{\ell}}=(Y_{\ell 1},\ldots,Y_{\ell N_{\ell}})^T$, we have

$$\mathrm{Var}(\varphi_{qn}(u,w))$$

$$\leq \frac{1}{n}E\left\{\frac{1}{E(N)}\sum_{j=1}^{N}\exp(iuT_j+iwZ)\psi_q(T_j,Z,Y_j)\right\}^2$$

$$\leq \frac{1}{n}E\left\{\left(\frac{1}{E(N)}\right)^2\left(\sum_{j=1}^{N}\exp(i2uT_j+i2wZ)\right)\left(\sum_{j=1}^{N}\psi_q^2(T_j,Z,Y_j)\right)\right\}$$

$$= \frac{1}{n}E(\psi_q^2(T,Z,Y)),$$

where the second inequality follows from Cauchy–Schwarz inequality. The lemma now follows from

$$E\left\{\sup_{t\in\mathcal{T};z\in\mathcal{Z}}|\Psi_{qn}-E\Psi_{qn}|\right\}$$

$$\leq \frac{\int\int E\{|\varphi_{qn}(u,w)-E(\varphi_{qn}(u,w))|\}|\zeta_1(h_{\mu,t}u,h_{\mu,z}w)|\,du\,dw}{4\pi^2 h_{\mu,t}^{\nu_1}h_{\mu,z}^{\nu_2}}$$

$$\leq \frac{1}{4\pi^2}\frac{\sqrt{E(\psi_q^2(T,Z,Y))}\int\int|\zeta_1(u,w)|\,du\,dw}{\sqrt{n}h_{\mu,t}^{\nu_1+1}h_{\mu,z}^{\nu_2+1}} = O\left(\frac{1}{\sqrt{n}h_{\mu}^{|\nu|+2}}\right). \quad\square$$

PROOF OF THEOREM 3.2. Let

$$S_{pq} = \sum_i\sum_j w_{ij}(T_{ij}-t)^p(Z_i-z)^q,$$

$$R_{pq} = \sum_i \sum_j w_{ij}(T_{ij} - t)^p (Z_i - z)^q Y_{ij},$$

where $w_{ij} = \frac{1}{nh_{\mu,t}h_{\mu,z}}K_2(\frac{t-T_{ij}}{h_{\mu,t}}, \frac{z-Z_i}{h_{\mu,z}})$. It can be shown that

$$\hat{\beta}_0 = \frac{R_{00} - \hat{\beta}_1 S_{10} - \hat{\beta}_2 S_{10}}{S_{00}},$$

where

$$\hat{\beta}_1 = \frac{-R_{00}(S_{10}S_{02} - S_{11}S_{01}) + R_{10}(S_{00}S_{02} - S_{01}^2) - R_{01}(S_{00}S_{11} - S_{01}S_{10})}{S_{00}(S_{10}S_{02} - S_{11}^2) - S_{10}(S_{10}S_{02} - S_{11}S_{01}) + S_{01}(S_{10}S_{11} - S_{01}S_{20})},$$

$$\hat{\beta}_2 = \frac{R_{00}(S_{10}S_{11} - S_{02}S_{20}) - R_{10}(S_{00}S_{11} - S_{10}S_{01}) + R_{01}(S_{00}S_{20} - S_{10}^2)}{S_{00}(S_{10}S_{02} - S_{11}^2) - S_{10}(S_{10}S_{02} - S_{11}S_{01}) + S_{01}(S_{10}S_{11} - S_{01}S_{20})}.$$

Applying Lemma C.2 and Slutsky's theorem repeatedly, we can show through tedious calculations and Theorem 3.1 that $|\hat{\beta}_1 - \beta_1| = O_p(\frac{1}{\sqrt{nE(N)h_\mu^4}})$ and $|\hat{\beta}_2 - \beta_2| = O_p(\frac{1}{\sqrt{nE(N)h_\mu^4}})$, where $h_\mu \asymp h_{\mu,t} \asymp h_{\mu,z}$. These results imply that

$$\lim \sqrt{n\bar{N}h_\mu^2}[(\hat{\mu}_L(t,z) - \mu(t,z)) - (\tilde{\mu}(t,z) - \mu(t,z))] \overset{\mathcal{D}}{\to} 0,$$

where $\tilde{\mu}(t,z) = [R_{00} - \beta_1 S_{10} - \beta_2 S_{10}]/S_{00}$. It suffices to show Theorem 3.2 for $\tilde{\mu}(t,z)$ instead of $\hat{\mu}_L(t,z)$, and this follows from setting $H(x_1, x_2, x_3, x_4) = [x_2 - \beta_1 x_3 + \beta_1 t x_1 - \beta_2 x_4 + \beta_2 z x_1]/x_1$, $\psi_1(u_1, u_2, u_3) = 1$, $\psi_2(u_1, u_2, u_3) = u_3$, $\psi_3(u_1, u_2, u_3) = u_1$, and $\psi_4(u_1, u_2, u_3) = u_2$ in Lemma C.2. □

## APPENDIX D: PROOFS OF THEOREMS 3.3 AND 3.4

For an integer $Q \geq 1$, let $\vartheta_q : R^5 \to R$ for $q = 1, \ldots, Q$ satisfy:

C.3 $\vartheta_q(t, s, z, y_1, y_2)$'s are continuous on $U(\{t, s, z\})$ uniformly in $(y_1, y_2) \in R^2$.

C.4 The functions $\frac{\partial^p}{\partial t^{p_1} \partial s^{p_2} \partial z^{p_3}} \vartheta_q(t, s, z, y_1, y_2)$ exist for all arguments $(t, s, z, y_1, y_2)$ and are continuous on $U(\{t, s, z\})$ uniformly in $(y_1, y_2) \in R^2$, for $p_1 + p_2 + p_3 = p$ and $0 \leq p_1, p_2, p_3 \leq p$.

The general weighted averages of three-dimensional smoothing methods are defined as

$$\Theta_{qn}(t, s, z) = \frac{1}{nE(N(N-1))h_{G,t}^{\nu_1 + \nu_2 + 2}h_{G,z}^{\nu_3 + 1}}$$

$$(D.1) \qquad \times \sum_{i=1}^n \sum_{1 \leq j \neq k \leq N_i} \vartheta_q(T_{ij}, T_{ik}, Z_i, Y_{ij}, Y_{ik})$$

$$\times K_3\left(\frac{t - T_{ij}}{h_{G,t}}, \frac{s - T_{ik}}{h_{G,t}}, \frac{z - Z_i}{h_{G,z}}\right),$$

where $K_3$ is a kernel function of order $(\boldsymbol{\nu}, \boldsymbol{\kappa})$ [see (A.2)]. Let

$$
\begin{aligned}
\xi_q(t,s,z) &= \frac{\partial^{|\boldsymbol{\nu}|}}{\partial t^{\nu_1} \, \partial s^{\nu_2} \, \partial z^{\nu_3}} \\
&\quad \times \int \vartheta_q(t,s,z,y_1,y_2) f_5(t,s,z,y_1,y_2) \, dy_1 \, dy_2, \\
\omega_{qr} &= \int \vartheta_q(t,s,z,y_1,y_2) \vartheta_r(t,s,z,y_1,y_2) \\
&\quad \times f_5(t,s,z,y_1,y_2) \, dy_1 \, dy_2 \|K_3\|^2,
\end{aligned}
$$

where $f_5(t,s,z,y_1,y_2)$ is the joint density of $(T_1, T_2, Z, Y_1, Y_2)$, $\|K_3\|^2 = \int K_3^2$, and $1 \leq q, r \leq l$.

Next, we provide the longitudinal version of the asymptotic normality property of kernel-weighted averages for three-dimensional smoothers.

LEMMA D.1. *Under Assumption* A.4–A.6, B.5–B.8, C.3 *and* C.4,

$$
\begin{aligned}
(\text{D.2}) \quad & \sqrt{nE[N(N-1)]h_{G,t}^{2\nu_1+2\nu_2+2} h_{G,z}^{2\nu_3+1}} \\
& \quad \times \{(\Theta_{1n}, \ldots, \Theta_{Qn})^T - (E(\Theta_{1n}), \ldots, E(\Theta_{Qn}))^T\} \\
& \xrightarrow{\mathcal{D}} N(0, \Omega).
\end{aligned}
$$

PROOF. The proof follows similar framework as in Lemma C.1 with appropriate modifications for three-dimensional smoothers. □

LEMMA D.2. *Let* $H : R^Q \to R$ *be a function with continuous first order derivatives,* $DH(v) = (\frac{\partial}{\partial x_1} H(v), \ldots, \frac{\partial}{\partial x_Q} H(v))^T$, *and* $\bar{N} = \frac{1}{n} \sum_{i=1}^n N_i$. *Under Assumption* A.4–A.6, B.5–B.8, C.3 *and* C.4, $\frac{h_{G,z}}{h_{G,t}} \to \rho_G$ *and* $nE(N(N-1))h_{G,t}^{2\boldsymbol{\kappa}+3} \to \tau_G^2$ *for some* $0 < \rho_G, \tau_G < \infty$, *we obtain*

$$
\begin{aligned}
& \sqrt{n\bar{N}(\bar{N}-1)h_{G,t}^{2\nu_1+2\nu_2+2} h_{G,z}^{2\nu_3+1}} \{H(\Theta_{1n}, \ldots, \Theta_{Qn}) - H(\xi_1, \ldots, \xi_Q)\} \\
& \xrightarrow{\mathcal{D}} N(\gamma_H, [DH(\xi_1, \ldots, \xi_Q)]^T \Omega [DH(\xi_1, \ldots, \xi_Q)]),
\end{aligned}
$$

*where* $\Omega = (\omega_{qr})_{1 \leq q, r \leq Q}$ *and*

$$
\begin{aligned}
\gamma_H = \sum_{q=1}^Q \sum_{\kappa_1 + \kappa_2 + \kappa_3 = \kappa} & \left\{ \frac{(-1)^{\boldsymbol{\kappa}}}{\boldsymbol{\kappa}!} \int u_1^{\kappa_1} u_2^{\kappa_2} u_3^{\kappa_3} K_3(u_1, u_2, u_3) \, du_1 \, du_2 \, du_3 \right\} \\
& \times \frac{d^{\boldsymbol{\kappa}}}{dt^{\kappa_1} \, ds^{\kappa_2} \, dz^{\kappa_3}} \int \vartheta_q(t,s,z,y_1,y_2)
\end{aligned}
$$

$$\times f_5(t, s, z, y_1, y_2) \, dy_1 \, dy_2$$

$$\times \frac{\partial H}{\partial \xi_q}(\xi_1, \ldots, \xi_Q)^T \tau_G \sqrt{\rho_G^{2\kappa_3+1}}.$$

PROOF. The framework of this proof is similar to that of Lemma D.1. □

LEMMA D.3. *Under assumptions* A.2, A.5–A.7, B.5–B.7, C.3 *and* C.4 *and* D.2,

$$\sup_{t,s \in \mathcal{T}; z \in \mathcal{Z}0} |\Theta_{qn} - \xi_q| = O_p\left(\frac{1}{\sqrt{n} h_G^{|\nu|+3}}\right) \qquad \text{where } h_G \asymp h_{G,t} \asymp h_{G,z}.$$

PROOF. The proof is very similar to the proof of Lemma C.3. □

LEMMA D.4. *Under assumptions* A.1 *and* A.2, A.5–A.7, B.1–B.3, B.5–B.7 *and* D.1 *and* D.2 *for* $K_2$, *we have defined as* (A.1) *of order* $(\mathbf{0}, 2)$, *and* $K_3$ *defined as* (A.2) *of order* $(\mathbf{0}, 2)$,

$$\sup_{t \in \mathcal{T}; z \in \mathcal{Z}} |\hat{\mu}_L(t, z) - \mu(t, z)| = O_p\left(\frac{1}{\sqrt{n} h_{\mu,t} h_{\mu,z}}\right),$$

$$\sup_{t,s \in \mathcal{T}; z \in \mathcal{Z}} |\hat{\Gamma}_L(t, s, z) - \Gamma(t, s, z)| = O_p\left(\frac{1}{\sqrt{n} h_{G,t}^2 h_{G,z}}\right).$$

PROOF. Apply Lemma C.3 to the Nadaraya–Watson estimator $\hat{\mu}_{\mathrm{NW}}(t, z)$ by choosing $\psi_1(t, z, y) = y$, $\psi_t(t, s, y) = 1$, and $H(x_1, x_2) = \frac{x_1}{x_2}$, one can obtain

$$\sup_{t \in \mathcal{T}; z \in \mathcal{Z}} |\hat{f}(t, z) - f(t, z)| = O_p\left(\frac{1}{\sqrt{n} h_\mu^2}\right),$$

$$\sup_{t \in \mathcal{T}; z \in \mathcal{Z}} |\hat{\mu}_{\mathrm{NW}}(t, z) - \mu(t, z)| = O_p\left(\frac{1}{\sqrt{n} h_\mu^2}\right).$$

Similar to the proof of Theorem 3.2, one can rewrite $\hat{\mu}_L(t, z)$ as

$$\hat{\mu}_L(t, z) = \hat{\mu}_{\mathrm{NW}}(t, z) - \frac{\hat{\beta}_1}{\hat{f}(t, z)} S_{10} - \frac{\hat{\beta}_2}{\hat{f}(t, z)} S_{10},$$

where $S_{10}$ is defined in the proof of Theorem 3.2 and show that

$$\sup_{t \in \mathcal{T}; z \in \mathcal{Z}} |\hat{\beta}_1 - \beta_1| = O_p\left(\frac{1}{\sqrt{n} h_\mu^3}\right) \quad \text{and} \quad \sup_{t \in \mathcal{T}; z \in \mathcal{Z}} |\hat{\beta}_2 - \beta_2| = O_p\left(\frac{1}{\sqrt{n} h_\mu^3}\right).$$

Thus, $\sup_{t \in \mathcal{T}; z \in \mathcal{Z}} |\hat{\mu}_L(t, z) - \mu(t, z)| = O_p(\frac{1}{\sqrt{n} h_\mu^3})$.

The uniform convergence rate of the covariance estimator simply replaces Lemma C.3 with Lemma D.3. □

PROOFS OF THEOREMS 3.3 AND 3.4. From Lemma D.4, we know that $\sup_{t,z}|\hat{\mu}(t,z) - \mu(t,z)| = O_p(\frac{1}{\sqrt{n}h_\mu^2})$ for both $\hat{\mu}_{\mathrm{NW}}(t,z)$ and $\hat{\mu}_L(t,z)$. Let $\vartheta_1(t_1,t_2,z,y_1,y_2) = (y_1 - \mu(t_1,z))(y_2 - \mu(t_2,z))$, $\vartheta_2(t_1,t_2,z,y_1,y_2) = (y_1 - \mu(t_1,z))$, and $\vartheta_3(t_1,t_2,z,y_1,y_2) = 1$, then $\sup_{t,z}|\Theta_{qn}| = O_p(1)$, $q = 1,2,3$, by Lemma D.3. Thus, we can obtain $\sup_{t,z}|\Theta_{qn}|O_p(\frac{1}{\sqrt{n}h_\mu^2}) = O_p(\frac{1}{\sqrt{n}h_\mu^2})$ for $q = 2,3$. By the fact that

$$
\begin{aligned}
C_{ijk} = {}&\tilde{C}_{ijk} + (Y_{ij} - \mu(T_{ij}, Z_i))(\mu(T_{ik}, Z_i) - \hat{\mu}(T_{ik}, Z_i)) \\
&+ (Y_{ik} - \mu(T_{ik}, Z_i))(\mu(T_{ij}, Z_i) - \hat{\mu}(T_{ij}, Z_i)) \\
&+ (\mu(T_{ij}, Z_i) - \hat{\mu}(T_{ij}, Z_i))(\mu(T_{ik}, Z_i) - \hat{\mu}(T_{ik}, Z_i)),
\end{aligned}
$$

and $\sup_{t,z}|\hat{\mu}(t,z) - \mu(t,z)|^2 = O_p(\frac{1}{nh_\mu^4})$ is negligible compared to $\Theta_{1n}$, $\hat{\Gamma}_{\mathrm{NW}}(t,s,z)$ and $\hat{\Gamma}_L(t,s,z)$ obtained via smoothing $C_{ijk}$ are asymptotically equivalent to those, denoted by $\tilde{\Gamma}_{\mathrm{NW}}(t,s,z)$ and $\tilde{\Gamma}_L(t,s,z)$, respectively, obtained via smoothing $\tilde{C}_{ijk}$. Therefore, it suffices to show the asymptotic distributions of $\tilde{\Gamma}_{\mathrm{NW}}(t,s,z)$ and $\tilde{\Gamma}_L(t,s,z)$.

Theorem 3.3 now follows from Lemma D.2 by letting $\vartheta_1(t,s,z,y_1,y_2) = (y_1 - \mu(t,z))(y_2 - \mu(s,z))$, $\vartheta_2(t,s,z,y_1,y_2) = 1$, and $H(x_1,x_2) = \frac{x_1}{x_2}$.

Theorem 3.4 follows from similar arguments as in the proof of Theorem 3.2. □

PROOF OF THEOREM 3.5. To show the asymptotic results of the mF-PCA covariance estimator, we need the following regularity conditions for the pooled covariance function and some joint p.d.f.'s:

E.1 $\frac{d^{\boldsymbol{\kappa}}}{dt^{k_1}ds^{k_2}} g_2(t,s)$ exists and is continuous on $\{(t,s)\}$ for $k_1 + k_2 = \boldsymbol{\kappa}$, $0 \le k_1, k_2 \le \boldsymbol{\kappa}$, and $g_2(t,s) > 0$;

E.2 $f_4(t,s,y_1,y_2)$ is continuous on $\{(t,s)\}$ uniformly in $(y_1,y_2) \in R^2$; $\frac{d^{\boldsymbol{\kappa}}}{dt^{k_1}ds^{k_2}} f_4(t,s,y_1,y_2)$ exists and is continuous on $\{(t,s)\}$ uniformly in $(y_1,y_2) \in R^2$, for $k_1 + k_2 = \boldsymbol{\kappa}$, $0 \le k_1, k_2 \le \boldsymbol{\kappa}$;

E.3 $f_8(t_1,t_2,t_1',t_2',y_1,y_2,y_1',y_2')$ is continuous on $\{(t_1,t_2,t_1',t_2')\}$ uniformly in $(y_1,y_2,y_1',y_2') \in R^4$.

Since the covariance estimator of mFPCA involves two-dimensional smoothing, the theoretical properties of covariance estimator in Theorem 1 in Yao, Müller and Wang (2005) and Theorem 2 in Yao (2007) can be applied directly. Thus, under assumptions A.5–A.7, D.1, E.1–E.3, $h_{G^*} \to 0$, $nh_{G^*}^6 \to \infty$ and $nh_{G^*}^8 < \infty$, one could obtain that

$$
\sup_{t,s \in \mathcal{T}} |\hat{\Gamma}^*(t,s) - \Gamma^*(t,s)| = O_p\left(\frac{1}{\sqrt{n}h_{G^*}^2}\right). \qquad \square
$$

**Acknowledgments.** The authors would like to express gratitude for the insightful comments of two referees and an Associate Editor.

## REFERENCES

Besse, P. and Ramsay, J. (1986). Principal components analysis of sampled functions. *Psychometrika* **51** 285–311. MR0848110

Bhattacharya, P. K. and Müller, H. G. (1993). Asymptotics for nonparametric regression. *Sankhyā Ser. A* **55** 420–441. MR1323398

Boente, G. and Fraiman, R. (2000). Kernel-based functional principal components. *Statist. Probab. Lett.* **48** 335–345. MR1771495

Bosq, D. (2000). *Linear Processes in Function Spaces: Theory and Application.* Springer, New York. MR1783138

Cardot, H. (2000). Nonparametric estimation of smoothed principal components analysis of sampled noisy functions. *J. Nonparametr. Stat.* **12** 503–538. MR1785396

Cardot, H. (2006). Conditional functional principal components analysis. *Scand. J. Statist.* **34** 317–335. MR2346642

Carey, J. R., Liedo, P., Müller, H. G., Wang, J. L., Sentürk, D. and Harshman, L. (2005). Biodemography of a long-lived tephritid: Reproduction and longevity in a large cohort of female mexican fruit flies, *Anastrepha Ludens. Experimental Gerontology* **40** 793–800.

Castro, P. E., Lawton, W. H. and Sylvestre, E. A. (1986). Principal modes of variation for processes with continuous sample curves. *Technometrics* **28** 329–337.

Chiou, J.-M., Müller, H.-G. and Wang, J.-L. (2003). Functional quasi-likelihood regression models with smooth random effects. *J. R. Stat. Soc. Ser. B Stat. Methodol.* **65** 405–423. MR1983755

Dauxois, J., Pousse, A. and Romain, Y. (1982). Asymptotic theory for the principal component analysis of a vector random function: Some applications of statistical inference. *J. Multivariate Anal.* **12** 136–154. MR0650934

Fan, J. and Gijbels, I. (1996). *Local Polynomial Modelling and Its Applications.* Chapman and Hall, London. MR1383587

Ferraty, F. and Vieu, P. (2006). *Nonparametric Functional Data Analysis: Theory and Practice.* Springer, New York. MR2229687

Guo, W. (2002). Funcitonal mixed effects models. *Biometrics* **58** 121–128. MR1891050

Hall, P. and Hosseini-Nasab, M. (2006). On properties of functional principal component analysis. *J. R. Stat. Soc. Ser. B Stat. Methodol.* **68** 109–126. MR2212577

Hall, P., Müller, H.-G. and Wang, J.-L. (2006). Properties of principal component methods for functional and longitudinal data analysis. *Ann. Statist.* **34** 1493–1517. MR2278365

James, G. M., Hastie, T. J. and Suger, C. A. (2000). Principal components models for sparse functional data. *Biometrika* **87** 587–602. MR1789811

Kneip, A. and Utikal, K. (2001). Inference for density families using functional principal component analysis. *J. Amer. Statist. Assoc.* **96** 519–532. MR1946423

Mas, A. and Menneteau, L. (2003). *High Dimensional Probability III.* 127–134. Birkhäuser, Basel. MR2033885

Paul, D. and Peng, J. (2009). Consistency of restricted maximum likelihood estimators of principal components. *Ann. Statist.* **37** 1229–1271. MR2509073

Peng, J. and Paul, D. (2009). A geometric approach to maximum likelihood estimation of the functional principal components from sparse longitudinal data. *J. Comput. Graph. Statist.* In press.

Ramsay, J. O. and Silverman, B. W. (2002). *Applied Functional Data Analysis: Methods and Case Studies.* Springer, New York. MR1910407

Ramsay, J. O. and Silverman, B. W. (2005). *Functional Data Analysis,* 2nd ed. Springer, New York. MR2168993

Rao, C. R. (1958). Some statistical methods for comparison of growth curves. *Biometrics* **14** 1–17.

Rice, J. and Silverman, B. (1991). Estimating the mean and covariance structure non-parametrically when the data are curves. *J. Roy. Statist. Soc. Ser. B* **53** 233–243. MR1094283

Rice, J. A. (2004). Functional and longitudinal data analysis: Prospectives on smoothing. *Statist. Sinica* **14** 631–647. MR2087966

Rice, J. A. and Wu, C. (2001). Nonparametric mixed effects models for unequally sampled noisy curves. *Biometrics* **57** 253–259. MR1833314

Shi, M., Weiss, R. and Taylor, J. (1996). An analysis of paediatric cd4 counts for acquired immune deficiency syndrome using flexible random curves. *J. Appl. Stat.* **45** 151–163.

Wang, Y. (1998). Mixed-effects smoothing spline anova. *J. R. Stat. Soc. Ser. C Stat. Methodol.* **60** 159–174. MR1625640

Wu, H. and Zhang, J.-T. (2006). *Nonparametric Regression Methods for Longitudinal Data Analysis: Mixed-Effects Modeling Approaches.* Wiley, Hoboken, NJ. MR2216899

Yao, F. (2007). Asymptotic distributions of nonparametric regression estimators for longitudinal of functional data. *J. Multivariate Anal.* **98** 40–56. MR2292916

Yao, F., Müller, H.-G. and Wang, J.-L. (2005). Functional data analysis for sparse longitudinal data. *J. Amer. Statist. Assoc.* **100** 577–590. MR2160561

Zhang, D., Lin, X., Raz, J. and Sowers, F. (1998). Semiparametric stochastic mixed models for longitudinal data. *J. Amer. Statist. Assoc.* **93** 710–719. MR1631369

Department of Statistics
University of California, Davis
Davis, California 95616
USA
E-mail: crjiang@wald.ucdavis.edu
           wang@wald.ucdavis.edu